\documentclass[11pt, a4paper]{article}

\usepackage{amsfonts,amsmath,amsthm,amssymb}
\usepackage{mathrsfs}
\usepackage[ruled, nothing]{algorithm}
\usepackage{epsfig}
\usepackage{enumerate}
\usepackage{a4wide}
\usepackage{umlaute}
\usepackage{dsfont}
\usepackage{mdwlist}
\usepackage{wasysym}

\usepackage{color}

\newtheorem{theorem}{Theorem}[section]

\newtheorem{lemma}[theorem]{Lemma}
\newtheorem{corollary}[theorem]{Corollary}
\newtheorem{definition}[theorem]{Definition}
\newtheorem{proposition}[theorem]{Proposition}

\newtheorem{conjecture}[theorem]{Conjecture}

\newcommand{\op}[1]{\text{\normalfont{#1}}}

\newcommand{\dist}[0]{\op{dist}}

\newcommand{\gap}[0]{\op{gap}}

\newcommand{\bB}[0]{\overline{b}}

\newcommand{\Prob}[1]{\op{Pr}\left[#1\right]}
\newcommand{\condProb}[2]{\op{Pr}\left[\left.#1 \; \right| \; #2 \right]}
\newcommand{\E}[1]{\mathbb{E}\left[#1\right]}

\newcommand{\distrib}[2]{\op{#1}\left(#2\right)}

\newcommand{\Bin}[1]{\distrib{Bin}{#1}}             

\newcommand{\bigO}[1]{\mathcal{O}\left(#1\right)}

\newcommand{\card}[1]{\left|#1\right|}

\newcommand{\io}{{\mathds 1}_\mathcal{O}}

\newcommand{\half}{\frac12}
\newcommand{\eps}{\varepsilon}

\title{Extremal Subgraphs of Random Graphs: an Extended Version\footnote{This is an extended version of our paper {\sl Extremal Subgraphs of Random Graphs}, submitted for publication, that contains more detailed remarks on the proof of Theorem~\ref{thm:cliques}.}}

\author{Graham Brightwell$^1$ \and Konstantinos Panagiotou$^{2}$ \and Angelika Steger$^3$ }

\begin{document}

\maketitle

\footnotetext[1]{Department of Mathematics, London School of Economics, Houghton Street, London WC2A 2AE, United Kingdom \texttt{G.R.Brightwell@lse.ac.uk}}
\footnotetext[2]{Max-Planck-Institute for Informatics, 66123 Saarbr\"ucken, Germany, \texttt{kpanagio@mpi-inf.mpg.de}}
\footnotetext[3]{Institute of Theoretical Computer Science, ETH Zurich, 8092 Zurich, Switzerland, \texttt{steger@inf.ethz.ch}}

\begin{abstract}
We prove that there is a constant $c >0$, such that whenever $p \ge n^{-c}$, with probability tending to 1 
when $n$ goes to infinity, every maximum triangle-free subgraph of the random graph $G_{n,p}$ is 
bipartite. This answers a question of Babai, Simonovits and Spencer~\cite{ar:bss90}. 

The proof is based on a tool of independent interest: we show, for instance, that the maximum cut of 
almost all graphs with $M$ edges, where $M\gg n$, is ``nearly unique''. More precisely, given a maximum 
cut $C$ of $G_{n,M}$, we can obtain all maximum cuts by moving at most $\mathcal{O}(\sqrt{n^3/M})$ 
vertices between the parts of $C$.
\end{abstract}
\numberwithin{equation}{section}

\section{Introduction}

It is well-known that, in many different contexts, large triangle-free 
graphs are bipartite. For example, Mantel~\cite{ar:m07} proved that the maximum triangle-free subgraph 
of a complete graph on~$n$ vertices is a complete bipartite graph with~$\lfloor n/2\rfloor$ vertices in 
one class and~$\lceil n/2\rceil$ vertices in the other class. 
Erd\H{o}s, Kleitman, and Rothschild~\cite{ar:ekr76} proved that such a statement is also true in a 
probabilistic sense. More precisely, they showed that, if~$T_n$ denotes a graph drawn uniformly at random 
from the set of all triangle-free graphs on~$n$ labeled vertices, then the probability that~$T_n$ is 
bipartite tends to 1 for~$n$ tending to infinity. This result was generalized 
independently by Steger~\cite{ar:s05} and Osthus, Pr\"omel and Taraz~\cite{ar:opt03} to the case that, 
in addition to the number of vertices, also the number of edges is prescribed. The following result is 
from~\cite{ar:opt03}.
\begin{theorem}\label{counting}
Let~$T_{n,m}$ denote a graph drawn uniformly at random from the set of all triangle-free graphs on~$n$ 
labeled vertices and~$m$ edges. Then, for any~$\varepsilon > 0$,
\[
\lim_{n\to\infty}\Prob{T_{n,m} \text{ is bipartite}}
=
\begin{cases}
	1,& \text{if } m = o(n) \\
	0,& \text{if } \frac{n}{2} \le m \le (1-\varepsilon)\frac{\sqrt{3}}{4}n^{\frac32}\sqrt{\log n} \\
	1,& \text{if } m \ge (1+\varepsilon)\frac{\sqrt{3}}{4}n^{\frac32}\sqrt{\log n}.
\end{cases}
\]
\end{theorem}

For a graph~$G$, let~$t(G)$ denote the maximum number of edges in a triangle-free subgraph (not necessarily 
induced) of~$G$, and let~$b(G)$ be the maximum number of edges in a bipartite subgraph of~$G$. So~$b(G)$ is 
just the maximum size of a {\em cut} in~$G$. Of course, we always have~$t(G) \ge b(G)$.  Our general 
intuition -- guided by the above results -- suggests that, for dense enough graphs, these two parameters 
will typically be equal.

In 1990, Babai, Simonovits and Spencer~\cite{ar:bss90} studied these parameters 
for the random graph~$G_{n,p}$, which was introduced by Erd{\H{o}}s and R\'enyi in \cite{ar:er60}. 
They proved, among others, the following result.  
\begin{theorem}\label{thm:T=B}
There is a positive constant~$\delta$ such that, for~$p \ge \half - \delta$, 
$$
\lim_{n \to\infty}\Prob{t(G_{n,p}) = b(G_{n,p})} = 1.
$$ 
\end{theorem}

It seems unlikely that the property ``$t(G_{n,p}) = b(G_{n,p})$'' has a threshold for constant~$p$; indeed, 
Babai et.\ al.\ asked in \cite{ar:bss90} whether this result could be extended to cover edge 
probabilities~$p$ of the form~$n^{-c}$, for some positive constant~$c$.

As far as we know, Theorem~\ref{thm:T=B} could hold whenever~$p=p(n) \ge n^{-1/2+\eps}$, for arbitrary 
\mbox{$\eps >0$}. The property does not hold, for example, when~$p_0(n) = \frac{1}{10} (\log n)^{1/2}n^{-1/2}$, 
as an easy calculation shows that the random graph~$G_{n,p_0}$ a.a.s.\ (asymptotically almost surely) has an 
induced 5-cycle~$H$ such that no other vertex has more than one neighbor in~$H$: any maximum-size 
triangle-free subgraph then includes all the edges of~$H$, and is not bipartite.    

In this paper, we answer affirmatively the question of Babai, Simonovits and Spencer: we prove that 
Theorem~\ref{thm:T=B} holds whenever~$p=p(n) \ge n^{-c}$, for some fixed~$c > 0$---our proof
gives $c=1/250$.  In fact, we prove the following stronger result.
\begin{theorem}\label{thm:extremalMain}
There is a positive constant~$c$ such that, if~$p=p(n) \ge n^{-c}$, then 
\[
\lim_{n\to\infty}\Prob{\text{every maximum triangle-free subgraph of~$G_{n,p}$ is bipartite}} = 1.
\]
\end{theorem} 

It should be noted that Theorem~\ref{counting} cannot be used directly to prove
a theorem of this type.  For given~$p$, the result does imply that there is an 
$m=m(n)$ such that the expected number of non-bipartite triangle-free subgraphs of 
$G_{n,p}$ with~$m$ edges is~$o(1)$, while the expected number of bipartite subgraphs 
with~$m$ edges tends to infinity.  However, the events that particular bipartite 
subgraphs exist in the graph are very far from being independent, so this certainly 
does not prove that there is a.a.s.\ a bipartite subgraph of~$G_{n,p}$ with this 
number~$m$ of edges.
                                                                                                   
Let us indicate our general strategy for proving Theorem~\ref{thm:extremalMain}, and explain the main 
points of difficulty. We need a little notation first. Let~$[n] := \{1,\dots,n\}$ and~$p = p(n) \ge n^{-c}$, 
where~$c>0$ is some fixed and small constant. For a bipartition~$\Pi = (A,B)$ of~$[n]$, and a graph~$G$ 
with vertex set~$[n]$, we let~$E(G;\,\Pi)$ denote the set of edges of~$G$ with one endpoint in each part.  
The edges of~$E(G;\, \Pi)$ are said to go {\em across}~$\Pi$; the other edges of~$G$ are said to be 
{\em inside} the (parts of the) partition.  A {\em~$d$-perturbation} of~$E(G;\, \Pi)$ is a triangle-free 
subgraph of~$G$ obtained from~$E(G;\, \Pi)$ by adding at most~$d$ edges inside~$\Pi$, and 
removing any number of edges.

An adaptation of the proof of Babai et. al. \cite{ar:bss90} enables us to restrict our attention to 
triangle-free subgraphs of the~$G_{n,p}$ that are ``almost bipartite'', specifically that 
are~$p^{-c}$-perturbations of some bipartite subgraph, for some positive constant~$c$. One of the new 
ingredients in this paper is the use of the sparse regularity lemma and a related probabilistic embedding 
lemma (see Section \ref{sec:preliminaries}) in order to cover cases where~$p=o(1)$.

If we now {\em fix} a partition~$\Pi$ with the additional constraint that the two classes~$A$ and~$B$ are 
roughly equal in size, it is not too hard to show that, with reasonably high probability, 
no~$p^{-c}$-perturbation of~$E(G_{n,p};\,\Pi)$ has more edges than~$E(G_{n,p};\,\Pi)$. However, for 
fixed~$G=G_{n,p}$, this is certainly not true simultaneously for \emph{all} partitions~$\Pi$: for instance, 
if~$\{x,y\}$ is an edge of~$G$, and~$x$,~$y$ and all their common neighbors are in~$A$, then~$E(G;\,\Pi)$ 
could be enlarged by adding the edge~$\{x, y\}$, keeping the graph triangle-free.  

On the other hand, we only need to consider partitions~$\Pi$ in which~$E(G;\,\Pi)$ is {\em optimal}, i.e., 
has the maximum number of edges among all bipartite subgraphs, or nearly so. By definition, such partitions 
have more edges going across them than typical partitions do, so it seems plausible that a fixed 
near-optimal~$E(G;\,\Pi)$ is still unlikely to have a~$p^{-c}$-perturbation with more edges. We are able to 
explicitly confirm this intuition.   

However, the calculations we are making will not work if there are too many near-optimal 
partitions~$\Pi$, as then it becomes too likely that {\em one} of them could be improved by 
a~$p^{-c}$-perturbation. The final ingredient of our proof is to show that this is unlikely to be the case: 
in the range we consider, a random graph typically has relatively few bipartitions that are optimal or 
near-optimal.

Before making this statement more precise, we need some more notation and terminology, which we shall
use throughout the paper. The \emph{distance} of two bipartitions/cuts~$\Pi = (A,B)$ and~$\Pi' = (A',B')$ 
of~$[n]$ is defined as the number of vertices in which they differ, i.e.
\begin{equation}\label{eq:defdist}
	\dist(\Pi, \Pi')
	:=
	\min
	\left\{
		|A' \cap A| + |B'\cap B|, |A'\cap B| + |B'\cap A|
	\right\}.
\end{equation}
Since
$
	(|A' \cap A| + |B'\cap B|) ~+~ (|A'\cap B| + |B'\cap A|) = n,
$
we have that~$\dist(\Pi, \Pi') \le \frac{n}2$ for all pairs $(\Pi,\Pi')$.  

An \emph{optimal} bipartition of a graph $G$ is a bipartition $\Pi$ such that the number $|E(G,\Pi)|$ 
of edges of $G$ across $\Pi$ attains the maximum value $b(G)$. 
It is convenient for us to consider the notion of a \emph{canonical} optimal bipartition $\Pi^*(G)$ 
of a graph $G$ on $[n]$: for instance, we fix a list of the bipartitions of $[n]$, and let 
$\Pi^*(G)$ be the first bipartition on the list that is optimal for $G$.  For any bipartition $\Pi$,
we then set~$\dist(G;\, \Pi) = \dist(\Pi,\Pi^*(G))$,  
the distance of~$\Pi$ to the canonical optimal bipartition of~$G$. 

We say that two cuts have \emph{gap}~$g$, if the difference of their sizes is precisely~$g$, i.e., if
\begin{equation}\label{eq:defgap}
	\gap(G;\, \Pi, \Pi') := \card{E(G;\, \Pi)} - |E(G;\, \Pi')| = g.
\end{equation}
Finally, we say that a cut~$\Pi$ has gap~$g$ if its number of edges differs from that of an optimal 
bipartition by exactly~$g$, so the gap of $\Pi$ is $\gap(G;\, \Pi) := b(G) - |E(G;\, \Pi)|$.  Our result 
for near-optimal bipartitions is as follows. We state it for the uniform random graph~$G_{n,M}$.

\begin{theorem}\label{thm:propertiesOfCuts}
There is a constant~$C > 1$ such that the following is true for sufficiently large~$n$. 
Let~$n^{-1} \ll p = p(n) \le \frac12$ and~$M = M(n) := p\binom{n}{2}$. Furthermore, let~$r = r(n)$ 
be a natural number satisfying~$r \ll (pn)^{1/8}$, let~$\omega = \omega(n)$ be any function tending to 
infinity with $n$, and set
\[
	s_0 := C \cdot \omega \cdot r^4 \cdot \sqrt{np^{-1}}.
\]
Then
\[
	\Prob{\exists \Pi : \gap(G_{n,M}; \, \Pi) = r-1 \;\text{ and }\; \dist(G_{n,M}; \, \Pi) \ge s_0}
	~\le~
	\omega^{-1}.
\]
\end{theorem}

In the case $r=1$, this result implies that, a.a.s., the distance between any two optimal bipartitions 
of $G_{n,M}$ is at most $2s_0$, where $s_0$ is any function tending to infinity faster than
$\sqrt{np^{-1}}$ (or $\sqrt{n^3/M}$).  In other words, for most graphs with $n$ vertices
and $M$ edges, the maximum cut is unique up to movements of a small number of vertices: any two 
maximum cuts have distance at most $\omega \sqrt{n^3/M}$.  It seems 
likely that in fact most graphs contains a ``backbone", i.e.\ a pair of ``large" sets $B,C$ that are 
contained in opposite parts in {\sl every} maximum cut. Note that this does not follow directly from 
Theorem~\ref{thm:propertiesOfCuts}.

Although our main focus is on the most appealing case of triangle-free graphs, our methods extend to more 
general settings. Let~$\mathcal{K}_\ell$ be the complete graph on~$\ell$ vertices. We have the following 
result, replacing the triangle by an arbitrary complete graph.
\begin{theorem}\label{thm:cliques}
Let~$\ell \ge 3$. There is a~$c = c(\ell) > 0$ such that, whenever~$p = p(n) \ge n^{-c}$,
$$
\lim_{n\to\infty}\Prob{\text{every maximum~$\mathcal{K}_\ell$-free subgraph of~$G_{n,p}$ is ($\ell$-1)-partite}} = 1.
$$ 
\end{theorem}
We believe that a similar result is true not only for complete graphs, but also for many other graphs as well. 
In a graph~$H$ with chromatic number~$\chi := \chi(H)$, a {\em color-critical edge} is an edge~$e$ such that 
the graph with edge set~$E(H) \setminus\{e\}$ has chromatic number~$\chi-1$. It is known that, if~$H$ has a 
color-critical edge, then the maximum number of edges in an~$H$-free graph is the Tur\'an number, i.e., the 
largest~$H$-free graph is the same as the largest~$\chi$-partite graph. If~$H$ does not have a color-critical 
edge, then this fails, as adding one edge to the Tur\'an graph does not create a copy of~$H$.

We expect that Theorem~\ref{thm:cliques} is true for any fixed~$H$ that has at least one color-critical edge. 
On the other hand, such a result automatically fails for any graph~$H$ without a color-critical edge. 
Babai et.\ al.\ \cite{ar:bss90} discuss what can be proved for graphs without color-critical edges. 
We treat neither case here. 

~\\
{\bf Outline of the Paper.} 
The paper is structured as follows. In Section \ref{sec:preliminaries}, we introduce some notation and 
state a few facts from the theory of random graphs. Let~$\mathcal{T}(G)$ denote the set of maximum 
triangle-free subgraphs of a given graph~$G$. In Section \ref{sec:nearOptimal}, we prove that, a.a.s., for 
every~$T\in\mathcal{T}(G_{n,p})$ there exists a bipartition~$\Pi_T$ such that~$T$ is 
a~$p^{-12}\log^2n$-perturbation of~$\Pi_T$. Next, in Sections~\ref{sec:cuts} and~\ref{sec:FKGproof}, we 
present the proofs of Theorems~\ref{thm:propertiesOfCuts} and \ref{thm:extremalMain}. 
We close in Section~\ref{sec:generalizations} with some remarks on the proof of Theorem~\ref{thm:cliques}.

\section{Preliminaries \& Notation}\label{sec:preliminaries}

In this section we will present some basic facts from the theory of random graphs and from probability theory, 
which we will use frequently in the remainder of the paper. 

Without further reference we will often use the following estimates for the tail of the binomial 
distribution, which can be found for instance in \cite{book:jlr00}.

\begin{lemma}\label{len:cherny}
Let~$X$ be a random variable that is binomially distributed with parameters~$n$ and~$p$, and 
set~$\lambda := \E{X} = np$. For any~$t\ge 0$, 
\[
	\Prob{X \ge \lambda + t} \le e^{-\frac{t^2}{2(\lambda + t/3)}}
	~~\text{ and }~~
	\Prob{X \le \lambda - t} \le e^{-\frac{t^2}{2\lambda}}.
\]
\end{lemma}

Let us introduce some additional notation. We denote by~$\mathcal{G}_n$ the set of all graphs with 
vertex set~$[n]$: our random graphs $G_{n,p}$ and $G_{n,M}$ will always be taken from $\mathcal{G}_n$. 
For~$G\in\mathcal{G}_n$ and~$X,Y \subseteq [n]$, we denote by~$E(G)$ the edge set of~$G$, by~$E(G;\,X)$ the 
set of edges between vertices in~$X$, and by~$E(G;\, X,Y)$ the set of edges of~$G$ joining a vertex of~$X$ 
and a vertex of~$Y$. Furthermore, we set~$e(G;\,X) := |E(G;\,X)|$ and~$e(G;\,X,Y) := |E(G;\, X,Y)|$, where 
edges inside~$X\cap Y$ are counted only once. 

Applying the above tail bounds to edge sets in random graphs, we easily obtain the following statement, which
is not best possible, but suffices for our purposes.  Unless stated otherwise, logarithms are always to 
the base~$e$.
\begin{proposition}\label{prop:folklore}
Let~$p \gg \frac{\log n}{n}$ and define
\[
\begin{split}
	\mathcal{B}_{n,p} := \Big\{ G&\in \mathcal{G}_n ~|~
																\exists X,Y \subseteq [n] \text{ such that } \\
																&X\cap Y = \emptyset, |X|\ge|Y|\ge 10p^{-1}\log n \text{ and }
																\big|e(G;X,Y) - p|X||Y|\big| \ge \frac12p|X||Y|
											\Big\}.
\end{split}
\]
Then
$
	\Prob{G_{n,p} \;\in\; \mathcal{B}_{n,p}} = o(1).
$
\end{proposition}

In the subsequent proofs, we will often exploit the ``equivalence'' of the binomial random graph 
model~$G_{n,p}$ and the uniform random graph model~$G_{n,M}$, when~$p = M/\binom{n}{2}$.  More precisely, 
we will use \emph{Pittel's inequality} (see e.g. \cite{book:jlr00}) which states that, for any 
property~$\mathcal{Q}$ of graphs,
\begin{equation}\label{eq:pittelsinequality}
	\Prob{G_{n,M} \not\in \mathcal{Q}} \le 3\sqrt{M}\Prob{G_{n,p} \not\in \mathcal{Q}}, ~~
	{\textstyle \text{ where~$p = M/\binom{n}{2}$}}.
\end{equation}

Now let us turn our attention to optimal bipartitions of the uniform random graph. Recall that, for a 
given graph~$G$, $b(G)$ denotes the number of edges in an optimal bipartition of~$G$; the following 
proposition provides bounds for~$b(G_{n,M})$, which hold with high probability.

\begin{proposition}\label{prop:mlcGnM}
Let~$M \gg n$. For sufficiently large~$n$,  
\[
\Prob{\frac{M}{2} \le b(G_{n,M}) \le \frac{M}{2} + \sqrt{4nM}} \ge 1 - e^{-n}.
\]
\end{proposition}

\begin{proof}
The inequality~$b(G_{n,M}) \ge M/2$ is well-known to hold for all graphs with~$M$ edges.  In order to 
show the upper bound, let~$p := M / \binom{n}{2}$,~$L := \frac{p}{2}\binom{n}{2}$ 
and~$\Delta := \sqrt{2pn^2(n-1)}$. In the sequel we will show
\begin{equation}\label{eq:mlcUpperBound}
	\Prob{b(G_{n,p}) \ge L + \Delta} \ll n^{-1}e^{-n},
\end{equation}
which, together with (\ref{eq:pittelsinequality}), proves the proposition. To see (\ref{eq:mlcUpperBound}), 
define for every partition~$\Pi = (A,B)$ of the vertex set of~$G_{n,p}$ the random variable 
\[
	X_\Pi = \begin{cases}
						1, & e(G_{n,p};\, \Pi) \ge L + \Delta \\
						0, & \text{otherwise.}
					\end{cases}
\]
\newcommand{\pe}[0]{\op{pe}}%
The number of edges of~$G_{n,p}$ across~$\Pi$ is binomially distributed with parameters~$|A||B|$ and~$p$; 
using Lemma \ref{len:cherny} we obtain, for sufficiently large~$n$,
\[
	\E{X_\Pi}
	=
	\Prob{\Bin{|A||B|, p} \ge L + \Delta}
	\le
	e^{-\frac{\Delta^2}{2(p|A||B| + \frac{\Delta}{3})}}
	\le
	e^{-\frac{\Delta^2}{2(L + \frac{\Delta}{3})}}
	\le
	e^{-2n}.
\]
Therefore, if we let~$X = \sum_\Pi X_\Pi$, we readily obtain~$\Prob{X = 0} \le 2^ne^{-2n} \ll {n^{-1}e^{-n}}$.
\end{proof}

Recall that a bipartition $\Pi$ has gap~$g$ if the number of the edges across $\Pi$ differs from the 
number of edges across an optimal partition by exactly~$g$.  The next proposition states that, a.a.s., all 
bipartitions of~$G_{n,M}$ with small gap are ``balanced''.

\begin{proposition}\label{prop:mlcIsBalanced}
Let~$M \gg n$,~$p := M / \binom{n}{2}$ and~$\lambda = \lambda(n) \ge 0$. Furthermore, let
\begin{eqnarray*}
	\mathcal{B}_{n,M}
	&:=&
	\Big\{G\in \mathcal{G}_{n,M} ~|~
		\text{every } \Pi = (A,B) \text{ with } \gap(G;\, \Pi) \le \lambda \text{ satisfies}
	\\
		&&\hspace*{1.4cm} \card{|A| - \frac{n}{2}} \le 3 n^{\frac34}p^{-\frac14} + \lambda^{\frac12}p^{-\frac12}
		\text{ and } \card{|B| - \frac{n}{2}} \le 3 n^{\frac34}p^{-\frac14} + \lambda^{\frac12}p^{-\frac12}
	\Big\}.
\end{eqnarray*}
For sufficiently large~$n$, we have 
$
	\Prob{G_{n,M} \in \mathcal{B}_{n,M}} \ge 1 - e^{-n}.
$
\end{proposition}
\begin{proof}
We show the analogous result for the binomial random graph~$G_{n,p}$ and use inequality (\ref{eq:pittelsinequality}) 
to prove the statement. Let~$\Pi = (A,B)$ be a partition of the vertex set and 
write~$|A| = \frac{n}{2} + d$ and~$|B| = \frac{n}{2} - d$. Now assume 
that~$|d| > 3n^{3/4}p^{-1/4} + \lambda^{1/2}p^{-1/2}$. The number of possible edges across~$\Pi$ is
\[
	|A| \cdot |B|
	=
	\left(\frac{n}{2} + d\right)\left(\frac{n}{2} - d\right)
	= 
	\frac{n^2}{4} - d^2
	\le
	\frac{n^2}{4} - (9n^{3/2}p^{-1/2} + \lambda p^{-1}).
\]
The number~$C_\Pi$ of edges across~$\Pi$ is binomially distributed with parameters~$|A||B|$ and~$p$. 
Let us assume that~$\Pi$ is a bipartition of~$G_{n,p}$ with gap at most~$\lambda$.  Using 
Lemma~\ref{len:cherny}, we obtain that, whenever~$n$ is sufficiently large, with probability larger 
than~$1 - e^{-\frac32n}$, the number of edges in~$G_{n,p}$ is at 
least~$p\binom{n}{2} - \sqrt{2}n^{3/2}p^{1/2}$.  Together with Proposition \ref{prop:mlcGnM}, this implies 
that every optimal bipartition of~$G_{n,p}$ contains, for sufficiently large~$n$, at 
least~$\frac{pn^2}{4} - n^{3/2}p^{1/2}$ edges. Hence, for sufficiently large~$n$, the probability 
that~$\Pi$ has gap less than~$\lambda$ is at most
\[
\begin{split}
	&\op{Pr}\Big[C_\Pi \ge \frac{pn^2}{4} - n^{3/2}p^{1/2} - \lambda\Big] \\
	&\hspace{-5mm}\le
	\op{Pr}\Big[C_\Pi \ge \E{C_\Pi} + (9n^{3/2}p^{-1/2} + \lambda p^{-1})\cdot p - (n^{3/2}p^{1/2} + \lambda)\Big] \\
	&\hspace{-5mm}\le
	\Prob{C_\Pi \ge \E{C_\Pi} + 8n^{3/2}p^{1/2}}
	\le e^{-2n},
\end{split}
\]
where the last step is again due to Lemma \ref{len:cherny}. Therefore,
\[
	\Prob{G_{n,p} \not\in \mathcal{B}_{n,p}}
	\le e^{-\frac32n} + 2^n\cdot e^{-2n}
	\ll n^{-1}e^{-n},
\]
and an application of Pittel's inequality (\ref{eq:pittelsinequality}) completes the proof.
\end{proof}

Finally, we state bounds for the number of \emph{non}-edges across any optimal bipartition of the random 
graph~$G_{n,M}$. The following corollary is a straightforward consequence of Propositions~\ref{prop:mlcGnM} 
and \ref{prop:mlcIsBalanced}.

\begin{corollary}\label{cor:bmlcGnM}
Let~$M \gg n$ and set
\[
	\bB(G) = \min\left\{|A||B| - e(G;\, \Pi) ~\big|~ \Pi = (A,B) \text{ is an optimal bipartition of }G\right\}.
\]
There is a constant~$C>0$ such that, for sufficiently large~$n$, 
\[
	\Prob{\bB(G_{n,M}) \ge \frac12\left(\binom{n}{2} - M\right) - \sqrt{\frac{Cn^5}{M}}}
	\ge
	1 - 2e^{-n}.
\]
\end{corollary}
\newcommand{\GnM}{{G_{n,M}}}
\newcommand{\GnMt}{{G_{n,M+t}}}
\begin{proof}
Set~$p := M / \binom{n}{2}$. First we apply Proposition \ref{prop:mlcGnM} to~$\GnM$ to obtain that, with 
probability at least~$1 - e^{-n}$, every maximum bipartition of~$\GnM$ has size smaller 
than~$\frac{M}{2} + \sqrt{4nM}$. Furthermore, we apply Proposition \ref{prop:mlcIsBalanced} with 
$\lambda = 0$ to obtain that, with probability larger than~$1 - e^{-n}$, all maximum cuts~$\Pi = (A,B)$ 
of~$\GnM$ satisfy~$|A|, |B| \ge \frac{n}{2} - 3n^{\frac34}p^{-\frac14}$. We deduce that, with probability 
at least~$1 - 2e^{-n}$, the minimum number of non-edges across any optimal bipartition is at least
\[
	\left(\frac{n}{2} - 3 n^{\frac34}p^{-\frac14}\right)\left(\frac{n}{2} + 3 n^{\frac34}p^{-\frac14}\right)
	- \frac{M}{2} - \sqrt{4nM},
\]
and the claim follows from~$p = M / \binom{n}{2}$ and~$M\le n^2$.
\end{proof}


\section{Finding a Near-Optimal Bipartition}\label{sec:nearOptimal}

\newcommand{\mF}{\mathcal{F}}
\newcommand{\mK}{\mathcal{K}}
Suppose we have~$p = p(n) \ge n^{-c}$ for some positive (small) constant~$c$. For a graph~$G$, we 
denote by~$\mathcal{T}(G)$ the set of maximum triangle-free subgraphs of~$G$. In this section, we will 
prove that, a.a.s., every~$T\in\mathcal{T}(G_{n,p})$ is ``almost'' bipartite. More precisely, our proof 
consists of two parts:
\begin{itemize}
\item In Lemma \ref{lem:pn2}, we mimic the proof of \cite{ar:bss90} to show that there is a 
bipartition~$\Pi = \Pi_T = (A,B)$ with at most~$o(pn^2)$ edges of~$T$ 
inside~$\Pi$, i.e., connecting vertices in~$A$ or in~$B$. The new ingredient here is an application of the 
sparse version of Szemer\'edi's regularity lemma and a probabilistic embedding lemma.
\item Second, in Lemma \ref{lem:almostlpartite}, we show that in fact there is a bipartition~$\Pi'$ with 
at most~$p^{-12}\log^2 n$ edges of~$T$ inside~$\Pi'$. This proof uses similar ideas as in \cite{ar:bss90}, 
but differs from the original proof in most details.
\end{itemize}
Before we continue with our proof, let us introduce a variant of 
\emph{Szem\'eredi's} \emph{regularity lemma} which can be meaningfully applied to sparse graphs. 
Before we state it formally, we need a few technical definitions.

\begin{definition}\label{def:eps_p_regular}
A bipartite graph ~$B = (V_1 \cup V_2, E)$  is called \emph{$(\eps,p)$-regular} if, for 
all~$V_1' \subseteq V_1$ and~$V_2' \subseteq V_2$  with ~$|V_1'| \ge \eps |V_1|$ 
and~$|V_2'| \ge \eps |V_2|$, 
\[
  \bigg| \frac{e(B;\,V_1', V_2')}{|V_1'||V_2'|} -\frac{|E|}{|V_1||V_2|} \bigg|
  \le \eps p.
\]
\end{definition}

\begin{definition}
Let~$G=(V,E)$ be a graph, and fix a constant~$\eps > 0$. A partition~$(C_i)_{i=0}^k$ of~$V$ is called an 
\emph{equitable partition with exceptional class~$C_0$} if~$|C_1|=|C_2|=\cdots=|C_k|$ 
and~$|C_0| \le \eps |C_1|$. An \emph{$(\eps,p)$-regular partition} is an equitable partition~$(C_i)_{i=0}^k$ 
such that, with the exception of at most~$\eps k^2$ pairs, the pairs~$(C_i,C_j)$, $1\leq i\leq j\leq k$, 
are~$(\eps,p)$-regular.
\end{definition}

\begin{definition}
Let~$G=(V,E)$ be a graph, and fix constants~$0<\eta\leq 1$,~$0<p\leq 1$ and~$b\geq 1$. We say that~$G$ is 
\emph{$(\eta,b,p)$-upper-uniform} if, for all disjoint sets~$X$ and~$Y$  with ~$|X|, |Y| \ge \eta |V|$,
\[
	\frac{e(G;\,X,Y)}{|X||Y|}\leq b p .
\]
\end{definition}
We now state the sparse variant of Szemer\'edi's regularity lemma; see \cite{inc:k97} and \cite{inc:kr03}.

\begin{theorem}\label{thm:szem_sparse}
For any~$0<\eps<1/2$ and~$b,m_0\geq 1$, there are constants~$\eta=\eta(\eps,b,m_0)>0$ 
and~$M_0=M_0(\eps,m_0)\geq m_0$ such that, for any~$p>0$, any~$(\eta,b,p)$-upper-uniform graph with at 
least~$m_0$ vertices has an~$(\eps,p)$-regular partition~$(C_i)_{i=0}^k$ such that~$m_0\leq k \leq M_0$.
\end{theorem}

A further tool which we will need in our proofs is an \emph{embedding lemma}, which essentially states 
that almost every graph that can be partitioned so that all pairs of classes are suitably dense 
and~$(\eps,p)$-regular contains a copy of any fixed graph~$H$. We need one further definition before we 
make this result precise.

\begin{definition}
For a graph~$H = (V_H, E_H)$ with vertex set~$V_H$ and edge set~$E_H$, let~$\mathcal{G}(H,n,m,\eps)$ be 
the class of graphs on vertex set~$V=\bigcup_{x\in V_H} V_x$, where the~$V_x$ are pairwise disjoint 
sets of size~$n$, and edge set~$E=\bigcup_{\{x,y\}\in E_H}E_{xy}$, where~$E_{xy}$ is the edge set of 
an~$(\eps, m/n^2)$-regular bipartite graph with~$m$ edges between~$V_x$ and~$V_y$. 
\end{definition}

Unfortunately, it can be shown that not all graphs in~$\mathcal{G}(H,n,m,\eps)$ contain a copy 
of~$H$. On the other hand, if~$m$ is sufficiently large and~$\eps$ is sufficiently small, we can hope 
that only a tiny fraction of the graphs in~$\mathcal{G}(H,n,m,\eps)$ do not contain a copy of~$H$. 
This was conjectured by Kohayakawa, {\L}uczak and R\"{o}dl in \cite{ar:klr97}.

\begin{conjecture}\label{conj:emb}
Let~$H$ be a fixed graph. For any~$\beta>0$, there exist constants~$\eps_0>0$,~$C>0$,~$n_0>0$ such that, 
for all~$m\geq Cn^{2-1/d_2(H)}$,~$n\geq n_0$, and~$0<\eps \leq\eps_0$, we have
\[
	\card{\{G\in \mathcal{G}(H,n,m,\eps): H \mbox{ is not a subgraph of } G\}}
	\le
	\beta^m\binom{n^2}{m}^{e(H)}.
\]
Here $d_2(H):=\max\left\{\frac{e_F-1}{|V(F)|-2}\mid F\subseteq H, |V(F)|\geq 3\right\}$ 
denotes the 2-density of a graph.
\end{conjecture}

In this work, we only need a weaker version of the above conjecture, where~$H$ is a complete 
graph and the number of edges~$m$ is slightly larger. The theorem below was proved by Gerke, 
Marciniszyn and Steger in \cite{inp:gms05}.

\begin{theorem}\label{thm:relaxedKLR}
Fix~$\ell \ge 3$. For all~$\beta > 0$, there exist constants~$n_0\in\mathbb{N}$,~$C > 0$, and~$\eps_0 > 0$ such that
\[
	\card{\{G\in \mathcal{G}(\mathcal{K}_\ell,n,m,\eps): \mathcal{K}_\ell \mbox{ is not a subgraph of } G\}}
	\le
	\beta^m \; \binom{n^2}{m}^{\binom{\ell}{2}},
\]
provided that~$m \ge Cn^{2 - 1/(\ell-1)}$,~$n\ge n_0$, and~$0 < \eps \le \eps_0$.
\end{theorem}

In fact, in \cite{inp:gms05} a much stronger \emph{counting} version of the above theorem was proved; we do not 
need this strengthening here. Note also that the above theorem implies Conjecture~\ref{conj:emb} 
for~$H = \mathcal{K}_3$, which was proved already by Kohayakawa, \L uczak and R\"odl in \cite{ar:klr96}.

A final ingredient in our proofs is the following lemma of Kohayakawa, R\"odl and Schacht \cite{ar:krs04}, which 
states that \mbox{$(\eps,p)$-regular} graphs, whose edge number is only specified within bounds, contain a 
\mbox{$(3\eps,p)$-regular} spanning subgraph with a given number of edges.

\begin{lemma}\label{lem:slicing}
Take~$p =p(n) \gg n^{-1}$. For every~$\eps > 0$,~$\alpha > 0$, and~$C > 1$, there exists an~$n_0$ such that the 
following holds. If~$B = (V_1\cup V_2, E)$ is an~$(\eps,p)$-regular graph satisfying~$|V_1|, |V_2| \ge n_0$ 
and~$\alpha p |V_1||V_2| \le e(B; V_1,V_2) \le C p |V_1||V_2|$, then there exists a~$(3\eps, p)$-regular 
graph~$B' = (V_1\cup V_2, E')$ with~$E'\subseteq E$ and~$|E| = \alpha p |V_1||V_2|$.
\end{lemma}

Now we proceed with our results. Recall that~$\mathcal{T}(G)$ denotes the set of maximum triangle-free 
subgraphs of the graph~$G$.

\begin{lemma}\label{lem:pn2}
Fix~$\eps > 0$.  There exists~$C > 0$ such that, for~$p\ge Cn^{-1/2}$, a random graph~$G_{n,p}$ a.a.s.\ has the 
following property.  For all~$T\in\mathcal{T}(G_{n,p})$ there is a partition~$\Pi_T = \Pi = (A,B)$ of the 
vertex set such that all but at most~$\eps pn^2$ edges of~$~T$ go across~$\Pi$. 
Furthermore,~$\frac{n}2 -\eps n\le |A|,|B| \le \frac{n}{2}+\eps n$.
\end{lemma}

\begin{proof}
The proof is similar to the proof of the analogous result in \cite{ar:bss90} for constant density~$p$. The 
new ingredients here are the sparse version of Szemer\'edi's regularity lemma (Theorem \ref{thm:szem_sparse}) 
and the probabilistic embedding lemma (Theorem \ref{thm:relaxedKLR}).

First we collect some properties of the random graph $G_{n,p}$. Using Chernoff's inequality, it is easy to 
verify that, for every~$c,\eps\in(0,1]$, a.a.s.\ every subset~$U$ of the vertices of~$G_{n,p}$ with~$|U| \ge cn$ 
spans more than~$(1-\eps) \frac12p|U|^2$ and less than~$(1+\eps) \frac12p|U|^2$ edges. Similarly, we have that, 
for every~$\xi,\eps>0$, a random graph~$G_{n,p}$ a.a.s.\ is such that, whenever~$X$ and~$Y$ are two disjoint 
subsets of the vertices with~$|X|,|Y| \ge \xi n$, we have~$|e(G_{n,p};\, X,Y) - p|X||Y|| \le \eps p|X||Y|$.
In particular, this implies that~$G_{n,p}$ is a.a.s.~$(\mu,(1+\eps),p)$-upper-uniform, for all fixed~$\mu>0$. 
Hence, a.a.s., Theorem~\ref{thm:szem_sparse} applies to~$G_{n,p}$ and all its spanning subgraphs.

Next we explain how to choose the constant~$C$. To do this we need some careful preparations.
Let 
\[
\mathcal{F}(n,m,\alpha) := \big\{G\in \mathcal{G}(\mathcal{K}_3,n,m,\alpha) ~:~ 
\mathcal{K}_3 \mbox{ is not a subgraph of } G\big\}.
\]
We apply Theorem \ref{thm:relaxedKLR} with~$\ell = 3$ and $\beta:=\frac{\eps^3}{e^6}$ to obtain the 
constants~$n_\eps$,~$C_\eps$ and~$\eps'$, which may depend on~$\eps$. Next we 
let~$\eps'' := \frac13\min\{\eps, \eps'\}$,~$b:= 1 + \eps$,~$m_0 := \eps^{-1}$, and apply 
Theorem~\ref{thm:szem_sparse} for~$\eps''$,~$b$ and~$m_0$ to obtain constants~$\eta$ and~$M_0$. 
Finally, we let~$\mu := \min\{\eta, \frac{1 - \eps''}{2M_0}\}$, and~$C := \frac{C_\eps}{\eps\mu^2}$.

We claim that, for all~$p\ge C n^{-1/2}$, the random graph~$G_{n,p}$ a.a.s.\ does not contain a graph 
from~$\bigcup_{\tilde n\ge\mu n}\mathcal{F}(\tilde{n},\eps p\tilde{n}^2,\eps')$. To see this, let~$X$ denote the 
number of such copies; we prove the claim by showing~$\E{X} = o(1)$. 
Let~$M(\tilde{n}) := \eps p \tilde{n}^2$ and observe that 
\begin{equation}\label{eq:numForbinGnpII}
	\mathbb E[X]
	\le
	\sum_{\tilde{n}\geq \mu n}
	n^{3\tilde{n}} \cdot |\mathcal{F}(\tilde{n},M(\tilde{n}),\eps')| \cdot p^{3M(\tilde{n})}.
\end{equation}
We assume now that $n \ge n_\eps/\mu$, so that each $\tilde n$ in the above sum is at least $n_\eps$.  
Now we recall that~$n_\eps, C_\eps$, and~$\eps'$ were chosen in such a way that we can apply 
Theorem~\ref{thm:relaxedKLR} with~$\beta:=\frac{\eps^3}{e^6}$ to obtain the 
bound~$|\mathcal{F}(\tilde{n},m,\eps')| \le \beta^{m}\binom{\tilde{n}^2}{m}^3$ whenever $\tilde n \ge n_\eps$ 
and~$m\ge C_\eps \tilde{n}^{3/2}$. We thus need to check that~$M(\tilde n)$ 
satisfies~$M(\tilde n)\ge C_\eps \tilde{n}^{3/2}$. This follows from our choice of~$C = \frac{C_\eps}{\eps\mu^2}$ 
and the assumption~$p\ge Cn^{-1/2}$: 
\[
	M(\tilde{n})
	= \eps p \tilde{n}^2
	\ge \eps p (\mu n)^2
	\ge \eps C \mu^2n^{3/2}
	\ge C_\eps n^{3/2}.
\]
Together with the inequality~$\binom{n}{k} \le (\frac{en}{k})^k$, we thus obtain from
(\ref{eq:numForbinGnpII}) that
\[
	\mathbb E[X]
	\le
	\sum_{\tilde{n}\geq \mu n}
		n^{3\tilde{n}}\beta^{M(\tilde{n})} \left(\frac{e}{\eps p}\right)^{3M(\tilde{n})} p^{3M(\tilde{n})}
	=
	\sum_{\tilde{n}\geq \mu n} n^{3\tilde{n}}	e^{-3M(\tilde{n})} =o(1),
\]
where the last two equalities follow from the choice of~$\beta$ and the fact that~$M(\tilde n) = \Omega(n^{3/2})$.  
This completes the proof of the claim.

Now consider a random graph~$G_{n,p}$ for~$p\ge Cn^{-1/2}$. The above discussion shows that~$G_{n,p}$ is 
a.a.s.~a $(\mu, 1+\eps,p)$-upper-uniform graph, and that it does not contain a graph from the 
set~$\mathcal{F}(\tilde{n},\eps p\tilde{n}^2,\eps')$, for all~$\tilde{n}\ge \mu n$. Furthermore, a.a.s., every 
subset~$U$ of the vertices of~$G_{n,p}$ with~$|U| \ge \mu n$ 
satisfies~$\big|e(G_{n,p};\, U) - \frac12 p|U|^2\big| \le \eps p |U|^2$, and, for every two disjoint subsets~$X$ 
and~$Y$ of size at least~$\mu n$, we have~$\big|e(G_{n,p};\, X,Y) - p|X||Y|\big| \le \eps p|X||Y|$.

In the remainder of the proof we assume that~$G_{n,p}$ has all these properties. Let~$T\in\mathcal{T}(G_{n,p})$ 
denote any maximum triangle-free subgraph of~$G_{n,p}$. We apply Theorem~\ref{thm:szem_sparse} to~$T$ 
with~$\eps''$,~$b = 1 + \eps$ and~$m_0 = \eps^{-1}$, to obtain an~$(\eps'',p)$-regular partition~$(C_i)_{i=0}^k$, 
where~$m_0 \le k \le M_0$.  Next we define the reduced graph~$R$ consisting of~$k$ labeled vertices corresponding 
to the classes~$C_1, \dots, C_k$, and an edge between two vertices whenever the corresponding partition classes 
form an~$(\eps'',p)$-regular bipartite graph with at least~$\eps p |C_1|^2$ edges. Now we show that, if~$R$ contains 
a triangle, then so does~$T$. To see this, observe first that we 
have~$|C_i| \ge \frac{(1-\eps'')n}{M_0} \ge \mu n$, for all~$1\le i\le k$. Additionally, if~$R$ contains a triangle, 
then by definition there exist three sets~$C_{i_1},C_{i_2},C_{i_3}$ inducing three bipartite graphs that 
are~$(\eps'',p)$-regular and contain at least~$\eps p |C_1|^2$ edges.  From Lemma~\ref{lem:slicing} we deduce 
that these bipartite graphs have (spanning) subgraphs with \emph{exactly}~$\eps p |C_1|^2$ edges, which 
are~$(3\eps'', p)$-regular and therefore also~$(\eps',p)$-regular. That is,~$T$ contains a graph 
from~$\mathcal{G}(\mathcal{K}_3,|C_1|,\eps p|C_1|^2,\eps')$.  As~$G_{n,p}$ and hence also~$T\subseteq G_{n,p}$ 
does not contain a graph from~$\mathcal{F}(|C_1|,\eps p|C_1|^2,\eps')$, this implies that~$T$ contains a 
triangle, contradicting the fact that~$T$ is triangle-free. We conclude that~$R$ contains no triangle.

The remainder of the proof is essentially the same as the proof of the \emph{Main Lemma} in \cite{ar:bss90} -- 
we only sketch roughly the details and refer the reader to \cite{ar:bss90} for a more detailed proof.

Since~$R$ contains no triangle, Tur\'an's theorem yields~$e(R) \le k^2/4$. On the other hand, we can 
show~$e(R) \ge (1 - 40\eps)\frac{k^2}{4}$. To see this, observe that the number of edges of~$T$ which join 
vertices of the same~$C_i$, or vertices of~$C_0$ to some other vertex, or correspond to a ``low-density'' or 
non-regular pair~$(C_i,C_j)$ is at most~$7\eps pn^2$. Furthermore, the number of edges of~$T$ in a 
``high-density'' regular pair~$(C_i,C_j)$ is at most~$(1+\eps)p\left(\frac{n}{k}\right)^2$. Therefore we obtain
\[
	e(T) \le e(R) \cdot p\left(\frac{n}{k}\right)^2 + 8\eps pn^2.
\]
But since our~$G_{n,p}$ a.a.s.\ has the property that any two disjoint sets~$X,Y$ of size~$n/2$ satisfy 
$e(G_{n,p};\,X,Y) \ge (1-\eps) p\frac{n^2}{4}$, we know that~$e(T) \ge (1-\eps)p \frac{n^2}{4}$,
from which the claimed lower bound for~$e(R)$ follows easily.  Now, due to the \emph{stability lemma} 
in~\cite{inc:s68}, there is a function~$\gamma\to 0$ (when~$\eps\to 0$), such that we can find a 
bipartition~$(A_R,B_R)$ of~$R$ with at most~$\gamma k^2$ edges inside the parts, 
and~$|A_R|, |B_R| \le \frac{k}{2} + \gamma k$. This completes the proof of the lemma, as it can easily be 
seen that this implies the existence of a bipartition of~$T$ with the claimed properties. 
(Note that the bound~$\pm \gamma k$ suffices to obtain the claim of the lemma if we start with a sufficiently 
small~$\eps>0$. We omit the details.)
\end{proof}

Before we proceed with showing that we can find a much better bipartition than the one guaranteed by the above 
lemma, we need two auxiliary tools, which will be used extensively in the sequel.

\begin{lemma}\label{lem:randomness}
Let~$k\ge 1$ be an integer,~$p \ge n^{-{1}/{3k}}$ and~$c\in(0,1)$. Then the random graph $G_{n,p}$ on the 
vertex set~$[n]$ has a.a.s.\ the following property.  For every subset~$U$ of the vertices of size~$\card{U} > cn$, 
there exists a set~$Q_U$ of~$\bigO{p^{-k}}$ vertices, such that every~$k$-tuple of~$[n] \setminus Q_U$ is 
completely joined to at least~$(1-c)p^k\card{U}$ and at most~$(1+c)p^k\card{U}$ vertices in~$U$.
\end{lemma}

\begin{proof}
The proof is very similar to the proof of the \emph{randomness lemma} in \cite{ar:bss90} and we omit some of the 
details. The important difference here is that the statement also holds for~$k$-tuples which might have a 
non-empty intersection with~$U$, and also that an upper bound on the size of the common neighborhood is given.

Fix a set $U$ of vertices with $|U| > cn$.
We call a set of~$k$ vertices~$\{v_1,\dots,v_k\}$ \emph{violating} if the number of their common neighbors 
in~$U$ is smaller than~$(1-c)p^k|U|$ or larger than~$(1+c)p^k|U|$.  Assume there exist~$t := Cp^{-k}$ pairwise 
disjoint sets of~$k$ vertices that are violating, where the constant~$C$ will be chosen later.  Let~$X$ denote 
the union of these sets, and let~$U' := U\setminus X$.  Observe 
that~$|X| = kt = Ckp^{-k} = \bigO{n^{1/3}} \ll p^k|U|=\Omega(n^{2/3})$.  For each of these~$t$ sets, the 
number of joint neighbors in~$U'$ is either smaller than~$(1-c)p^k|U| \le (1-c/2)p^k|U'|$ or larger 
than~$(1+c)p^k|U| - |X| \ge (1+c/2)p^k|U'|$. For each~$k$-tuple, the probability that the number of joint 
neighbors in~$U'$ is so small or so big is, by Chernoff's theorem, less than~$e^{-C'p^kn}$, for an appropriate 
constant~$C'$. As all these events are independent, we obtain that the probability that there exists a set $U$ 
with~$t$ violating sets can be bounded from above by~$2^n$ (the number of ways to choose~$U$) times~$n^{kt}$ 
(the number of ways to choose the~$t$ sets of~$k$ vertices each) times~$e^{-C'tp^kn}$; hence the probability 
that~$G_{n,p}$ does not satisfy the claim of the lemma is at most
\begin{equation}\label{eq:l3.9}
	2^n \cdot n^{kt} \cdot e^{-C'tp^kn} = o(1),
\end{equation}
if~$C$ is chosen appropriately.
\end{proof}

The above lemma captures only cases where the set~$U$ has linear size. For sets $U$
of size $o(n)$, the above proof works similarly as before -- we just have to replace the left hand side of 
$(\ref{eq:l3.9})$ by $n^{|U|} \cdot n^{kt} \cdot e^{-C'tp^k|U|}$, from which we see that, for sets~$U$ of 
sublinear size, we have the following lemma. 

\begin{lemma}\label{lem:randomnessGraham}
Fix~$c\in (0,1)$ and let $k\ge 1$ be an integer. Then there is a~$C = C(c,k)>0$ such that~$G_{n,p}$ has 
a.a.s.\ the following property. For every subset~$U$ of the vertices with~$p^k|U| > C \log n$, there exists 
a set~$Q_U$ of at most $Cp^{-k}\log n$ vertices, such that every~$k$-tuple of~$[n] \setminus Q_U$ is 
completely joined to at least~$(1-c)p^k\card{U}$ and at most~$(1+c)p^k\card{U}$ vertices in~$U$.
\end{lemma}

With these preparations we can now prove the main result of this section.

\begin{lemma}\label{lem:almostlpartite}
For~$p\ge n^{-1/7}$, the following holds a.a.s.: for every~$T\in\mathcal{T}(G_{n,p})$, there is a 
partition~$\Pi_T= \Pi = (A,B)$ of the vertex set with at most~$p^{-12}\log^2 n$ edges of~$~T$ inside~$\Pi$. 
Furthermore,~$|A| = \frac{n}{2} + o(n)$ and~$|B| = \frac{n}{2} + o(n)$.
\end{lemma}

\begin{proof}
\newcommand{\GGnp}{\Gamma_{G_{n,p}}}
\newcommand{\Gnp}{{G_{n,p}}}
We restrict our proof to the case~$p = o(1)$, because for the remaining cases the statement follows directly 
from results in \cite{ar:bss90}.  Note that the assumption~$p=o(1)$ implies that, for~$n$ sufficiently large, 
$p^{-1}$ is larger than any given constant. Within the proof we will often use this fact in order to keep the 
formulas simpler.  A more careful handling of the inequalities would easily lead to an improvement on the 
bound~$p\ge n^{-1/15}$ in the statement of the lemma. As however our bound in the second part of the proof 
(Sections~\ref{sec:cuts} and~\ref{sec:FKGproof}) is even weaker, we put emphasis on readability of the proof 
instead of optimizing the constant.  

Let us first collect some properties holding a.a.s.\ for a random graph~$\Gnp$ with vertex set $[n]$.  
We apply Lemma~\ref{lem:randomness} with~$k=2$ and~$c=\frac14$, which yields that~$G_{n,p}$ has a.a.s.\ the 
following property:
\begin{itemize*}
\item[(P1)] For every set~$U$ of size at least~$\frac{n}4$, all pairs of vertices have between~$\frac34p^2|U|$ 
and~$\frac54p^2|U|$ common neighbors in~$U$, except for pairs intersecting a set $Q_U$ of~$\bigO{p^{-2}}$ 
vertices.
\end{itemize*}
Next, we apply Lemma~\ref{lem:randomnessGraham} with~$c= \frac14$ and $k=1$, obtaining 
that, for a suitable constant $C$, $\Gnp$ a.a.s.\ satisfies the following: 
\begin{itemize*}
\item[(P2)] For every subset~$U$ of the vertex set with $|U| > Cp^{-1}\log n$, 
the number of vertices having more than~$\frac54p|U|$ or less than~$\frac34p|U|$ neighbors in~$U$ is at 
most~$Cp^{-1}\log n$.
\end{itemize*}
We apply Proposition~\ref{prop:folklore} and obtain that~$\Gnp$ has a.a.s.\ the following property: 
\begin{itemize*}
\item[(P3)] For any two disjoint subsets $X$ and $Y$ of the vertex set, 
each of size at least $10p^{-1}\log n$, the number of edges between $X$ and $Y$ lies between 
$\frac12 p |X|\,|Y|$ and $\frac32 p |X|\,|Y|$.
\end{itemize*}

We now fix a constant $\eps$ with~$0 < \eps < \frac{1}{1000}$, and apply Lemma~\ref{lem:pn2} using~$\eps^5$ 
in place of~$\eps$, to obtain that $G_{n,p}$ satisfies the following property a.a.s.: 
\begin{itemize*}
\item[(P4)] For all~$T\in\mathcal{T}(G_{n,p})$, there is a partition~$\Pi_T = (A_T,B_T)$ 
such that $e(T;\,A_T) + e(T;\,B_T) \le {\eps^5} pn^2$ 
and~$\frac{n}{2} - \eps^5 n \le |A_T|,|B_T| \le \frac{n}{2} + \eps^5 n$.  
\end{itemize*}
Finally, we apply Proposition~\ref{prop:mlcIsBalanced} with~$\lambda = \eps^5pn^2$.  This tells us that, 
a.a.s., every bipartition with gap at most $\eps^5 pn^2$ has parts of size at least 
$\frac n2 - n^{3/4}p^{-1/4} - \lambda^{1/2}p^{-1/2}$.  For the given value of $\lambda$, this yields:  
\begin{itemize*}
\item [(P5)] Every bipartition~$\Pi' = (A',B')$ with gap at most~$\eps^5 pn^2$ 
satisfies~$(1 - \eps)\frac{n}{2} \le |A'|,|B'| \le (1 + \eps)\frac{n}{2}$.
\end{itemize*}

In the remainder of the proof, we show that, for $n$ sufficiently large, any graph~$G$ on $[n]$ satisfying all of 
(P1)-(P5), for the given values of $p$, $\eps$ and $C$, satisfies the conclusion of the lemma.
This will clearly suffice to complete the proof.  

Fix then such a graph $G$, and take any 
$T\in\mathcal{T}(G)$.  We shall call a partition~$\Pi = (A,B)$ \emph{optimal with respect to~$T$} 
if~$e(T; \, \Pi)$ attains its maximum over all possible partitions. Recall that~$b(G)$ denotes the size of an 
optimal bipartition of~$G$.  From (P4) we know that, if~$\Pi$ is optimal with respect to~$T$, we have
\[
	e(T;\,\Pi)
	\ge e(T) - \eps^5pn^2
	\ge b(G) - \eps^5pn^2,
\]
which implies that~$\Pi$ has gap at most~$\eps^5pn^2$.  Using (P5), we deduce that all optimal 
bipartitions~$\Pi = (A,B)$ of~$T$ have the properties~$(1-\eps)\frac{n}{2} \le |A|,|B| \le (1+\eps)\frac{n}{2}$ 
and~$e(T;\, A) + e(T;\,B) \le {\eps^5} pn^2$.  From now on, we fix such an optimal 
bipartition~$\Pi = (A,B)$ of $T$.
 
Before we continue, let us introduce some notation. For a graph~$J$, and a vertex~$v$, let $\Gamma(J;v)$ 
denote the set of neighbors of~$v$ in~$J$ and, for $S$ a set of vertices of $J$, define 
\[
	d(J;\, v,S) := \card{\Gamma(J;v) \cap S}.  
\]
We call an edge \emph{horizontal} if it is in~$T$ and joins two vertices in~$A$ or two vertices in~$B$. 
The \emph{horizontal degree} of a vertex~$v\in A$ is given by
\[
	d_H(v) := d(T; v, A).
\]
We call an edge \emph{missing} if it is in~$G$, joining two vertices in~$A$ and~$B$, \emph{but is not} in~$T$. 
The number of missing edges at a vertex~$v\in A$ (with respect to the partition~$\Pi$) is thus
\[
	d_M(v) := d(G; v, B) - d(T; v, B).
\]
In the remainder of the proof we will repeatedly use the following strategy. We will assume that the horizontal 
edges satisfy some property.  We then will use the assumptions (P1)-(P3) on~$G$, and the assumption that~$\Pi$ 
is optimal with respect to~$T$, in order to expose more missing edges than horizontal edges. However, this will 
clearly contradict the maximality of~$T$, as we could delete all horizontal edges from~$T$, and add all missing 
edges to it in order to obtain a larger triangle-free graph.

In order to formalize the idea, let us first define some sets of \emph{exceptional} vertices and discuss a 
few of their properties.  For a subset~$U$ of the vertex set $[n]$ of $G$ that contains at least $\frac{n}4$ vertices we set 
\[
\mathcal{B}_1(U) := \left\{
v \in [n] : \big||\Gamma(G; v) \cap U| - p|U|\big| \ge \frac{p|U|}4
\right\}.
\]
and
we set ${\mathcal B}_2(U)$ equal to the set $Q_U$ guaranteed by property (P1), so that any two elements
of $[n]\setminus {\mathcal B}_2(U)$ have between $\frac34 p^2 |U|$ and $\frac54 p^2 |U|$ common 
neighbors in $U$.  
Now set 
\begin{eqnarray}
	X_1^A
	&=&
	\big(\mathcal{B}_1(A) \cup \mathcal{B}_1(B) \cup \mathcal{B}_2(A) \cup \mathcal{B}_2(B)\big) \;\cap\; A,
	\label{eq:exceptional:randomness} \\
	X_2^A
	&=&
	\big\{
		v\in A ~|~ d_H(v) \ge \eps pn
	\big\} \setminus X_1^A,
	\label{eq:exceptional:largeDeg} \\
	X_3^{A}
	&=&
	\big\{
		v\in A ~|~ d_M(v) \ge d_H(v) + 5\eps pn
	\big\}
	\setminus (X_1^A\cup X_2^A).
	\label{eq:exceptional:missing}
\end{eqnarray}
Further, set~$A_0 = A$ and~$A_i = A_{i-1} \setminus X_i^A$ for~$i\in\{1,2,3\}$.  We also make the analogous 
definitions for each of the above sets with~$A$ replaced by~$B$.  Finally, let~$X_i = X_i^A\cup X_i^B$ 
for~$i\in\{1,2,3\}$.

Observe that, from (P1) and (P2), we have~$|X_1| = \bigO{p^{-2}}$.  Moreover, we can estimate the number of 
vertices in~$X_3$ as follows.  The number~$m$ of missing edges in~$T$ incident to at least one vertex 
in~$X_3$ satisfies
\[
	m \ge \frac12 \sum_{v\in X_3} d_M(v) \ge \frac12 |X_3| \cdot 5\eps pn.
\]
But~$m$ is no greater than~$\eps^5pn^2$
and we deduce~$\card{X_3} \le \eps n$, with room to spare.

We now deduce the following statements in turn, which gather ``self-improving'' information on the number of 
horizontal edges.
\begin{itemize*}
	\item[(i)] $X_2$ is small, i.e.,~$|X_2| \le \eps p^{-2}$.
	\item[(ii)] Set~$H_3 := E(T;\, A_3) \cup E(T;\, B_3)$. For large~$n$ we 
		have~$\card{H_3} \le p^{-2}n\log n$.
	\item[(iii)] $\card{X_3} \le p^{-4}\log n$ (i.e.\ the actual number of exceptional vertices 
		in~$X_3$ is much smaller than is guaranteed by our earlier calculation).
	\item[(iv)] There are no vertices in~$A_3$ or~$B_3$ with horizontal degree in~$H_3$ greater 
		than~$p^{-6}\log n$.
	\item[(v)] We improve the bound on $|H_3|$:~$|H_3| \le \frac12p^{-12}\log^2n$.
	\item[(vi)] Finally, $e(T; A) + e(T; B) \le p^{-12}\log^2 n$.
\end{itemize*}

\medskip
\noindent {\bf (i)} \quad 
Let~$X_2^A = \{v_1, \dots, v_t\}$, where~$t = |X_2^A|$, and observe that for every~$v\in X_2^A$ we 
have~$d_H(v) \le d(T; v, B)$, as otherwise~$\Pi$ would not be optimal -- we could improve $\Pi$ by moving 
$v$ to $B$. 
So $|\Gamma(T;v) \cap A|$ and $|\Gamma(T;v) \cap B|$ are both at least $\eps pn$, for each $v \in X_2^A$.  
Furthermore, all edges between the sets~$\Gamma(T;v)\cap A$ and~$\Gamma(T; v)\cap B$ in~$G$ are missing, 
as otherwise~$T$ would contain a triangle. 

In the following we show a lower bound for the total number of missing edges~$m$, which will immediately 
translate into an upper bound for~$t$. We write
\[
	\mathcal{F}(v_i) := E\big(G;~\Gamma(T; v_i)\cap A,\;\Gamma(T; v_i)\cap B\big),
\]
i.e.,~$\mathcal{F}(v_i)$ is the set of edges in~$G_{n,p}$ that are missing ``due to''~$v_i$.  We 
have~$f_i := |\mathcal{F}(v_i)| \ge \frac{p}{2} \cdot (\eps pn)^2$, from (P3). 
Let~$t_0 = (\frac{\eps}{2p})^2$, and suppose that~$t$ 
satisfies~$t \ge t_0$. Now set~$m_0 := |\bigcup_{i=1}^{t_0} \mathcal{F}(v_i)|$, and observe that~$m \ge m_0$. 
In order to bound~$m_0$, we apply the inclusion-exclusion principle:
\begin{equation}\label{eq:mr}
	m_0
	\ge
	\sum_{i=1}^{t_0} f_i - \sum_{1\le i < j\le t_0} \card{\mathcal{F}(v_i) \cap \mathcal{F}(v_j)}
	\ge
	t_0 \cdot \frac{p}{2} \cdot (\eps pn)^2 - \binom{t_0}{2}\max_{i < j}\card{\mathcal{F}(v_i) \cap \mathcal{F}(v_j)}.
\end{equation}
The size of the common neighborhood of any two vertices~$v,w\in X_2^A$ in~$A$, as well as in~$B$, is at 
most~$p^2n$, by the definition of~$X_1$.  So~$\card{\mathcal{F}(v) \cap \mathcal{F}(w)}$ is at most the 
maximum number of edges of~$G$ between two disjoint vertex sets of this size.  Hence, again by (P3), we have
\begin{equation}\label{eq:zzz}
	\card{\mathcal{F}(v) \cap \mathcal{F}(w)}
	\le 2p \cdot (p^2n)^2.
\end{equation}
 From (\ref{eq:mr}), (\ref{eq:zzz}) and the definition of $t_0$, we obtain~$m_0 \ge \frac{\eps^4}{16}pn^2$. 
Therefore, whenever~$|X_2^A| \ge t_0$, we achieve a contradiction, as~$m_0$ is at most~$\eps^5pn^2$.  
Similarly $|X_2^B| < t_0$, so we have $|X_2| < 2t_0 \le \eps p^{-2}$, as claimed.  

\medskip
\noindent {\bf (ii)} \quad 
Recall that~$A_3$ and~$B_3$ denote the sets of vertices which are not exceptional. 
We may assume~$e(T;\, A_3) \ge e(T;\, B_3)$, as otherwise we could interchange the roles of~$A$ and~$B$. Our 
objective is to derive upper and lower bounds for the number~$N$ of instances of a configuration called 
a ``chord''; these bounds will immediately imply a bound on~$|H_3|$ that will show (ii).  A chord consists 
of three vertices~$x,y\in A_3$ and~$z\in B$ with the property that~$x$ and~$y$ are connected by an edge 
in~$H_3$,~$y$ and~$z$ are connected in~$G$, and the edge~$\{x,z\}$ is missing (i.e., it is in~$G$ but not 
in~$T$).

\begin{figure}[ht]
\label{fig:3chord}
\centering
\includegraphics[height=34mm]{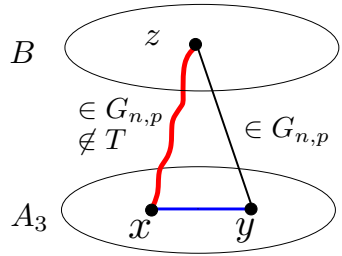}
\caption{A chord in $T$.}
\end{figure}

Consider an edge~$\{x,y\} \in E(T;\,A_3)$. Then every vertex~$z$ in the common neighborhood of~$x$ and~$y$ 
in~$B$ forms a triangle in~$G$.  Hence at least one of the edges~$\{x,z\}$ or~$\{y,z\}$ is a missing edge 
(as otherwise there would be a triangle in $T$), which means that the three vertices~$x$,~$y$ and~$z$ form 
a chord.  Due to the definition of~$X_1$, the common neighborhood of any two vertices of~$A_3$ in~$B$ has 
size at least~$\frac34p^2|B|$.  Therefore, a lower bound for the number of chords is
\begin{equation} \label{eq:NlowerBound}
	N
	\;\ge\;
	e(T; A_3) \cdot \frac34p^2|B|
	\;\ge\; \frac12|H_3|\cdot \frac{p^2n}{4}
	\;=\; |H_3|\frac{p^2n}{8}.
\end{equation}
In the sequel we derive an upper bound for~$N$. We first bound the number of chords that contain a 
vertex~$x\in A_3$ such that~$d_{H_3}(x) \le C p^{-1}\log n$, where~$C$ is the constant in (P2).  We bound 
the number of such chords containing~$x$ by their horizontal degree times an upper bound on the size of a 
common neighborhood of any two vertices in~$A_3$ in~$G$.  Using again the definition of~$X_1^A$, 
we deduce that the common neighborhood of any two vertices in~$A_3$ in~$B$ is of size at 
most~$\frac54p^2|B|$. Hence,
\begin{equation}\label{eq:NupperBoundFirst}
\begin{split}
	&\sum_{\substack{x \in A_3: \\ d_{H_3}(x) \le Cp^{-1}\log n}}
			\sum_{w \in \Gamma(T;\, x) \cap A_3} |\Gamma(G;\, x) \cap \Gamma(G;\, w) \cap B| \\
	&~~~~~~~~\le
	|A_3| \cdot Cp^{-1}\log n \cdot \frac54p^2|B|
	\;\le\; Cp\;n^2\log n.
\end{split}
\end{equation}
Now consider the vertices~$x\in A_3$ such that~$d_{H_3}(x) \ge Cp^{-1}\log n$.  The number of chords 
containing such a vertex can be bounded as follows. Consider the set of missing neighbors of~$x$, 
i.e.,~$S := (\Gamma(G;\, x) \cap B)\setminus \Gamma(T;\, x)$, and note that~$|S| = d_M(x)$.  By (P2), all 
but at most~$Cp^{-1}\log n$ vertices in~$S$ have in~$G$ at most~$\frac54pd_{H_3}(x)$ neighbors in the 
set~$\Gamma(T;\, x) \cap A_3$.  Furthermore, note that~$d_M(x) \le d_{H_3}(x) + 5\eps pn \le 6\eps pn$, 
as~$x\in A_3$.  We deduce that the number of chords containing a vertex of high horizontal degree can be 
bounded by
\begin{equation}\label{eq:NupperBoundSec}
	\begin{split}
		&\sum_{x: d_{H_3}(x) \ge Cp^{-1}\log n} \left(d_{M}(x)\cdot \frac54pd_{H_3}(x) + 
		Cp^{-1}\log n \cdot d_{H_3}(x)\right) \\
		&\le
		\sum_{\substack{x \in A_3: \\ d_{H_3}(x) \ge Cp^{-1}\log n}}
			\Big(\frac{30}{4}\eps p^2 n + Cp^{-1}\log n\Big) \cdot d_{H_3}(x)
		\le 16\eps p^2n \cdot |H_3|.
	\end{split}
\end{equation}	
Now by combining (\ref{eq:NlowerBound}), (\ref{eq:NupperBoundFirst}) and (\ref{eq:NupperBoundSec}), 
we obtain, for sufficiently large~$n$,
\[
\frac{p^2 n}8 |H_3| \le N \le Cpn^2 \log n + 16 \eps p^2 n |H_3|,
\]
and so 
\[
	|H_3| \cdot p^2n\left(\frac18 - 16\eps\right)
	\le Cp\;n^2\log n,
\]
from which the bound on~$|H_3|$ follows, with room to spare, from our assumptions on~$p$ and~$\eps$.

\medskip
\noindent {\bf (iii)} \quad 
Let~$H = E(T;\, A) \cup E(T;\, B)$, and recall that~$H_3$ is the subset of~$H$ restricted 
to~$A_3$ and~$B_3$. Our strategy is as follows.  We estimate~$|H|$ from above and the number of missing 
edges from below.  Comparing the two bounds will yield a contradiction, unless~$\card{X_3} \le p^{-4}\log n$. 
First observe that, due to (i) and (ii), the number of horizontal edges is at most
\begin{equation}\label{eq:exploitMissingI}
	\sum_{v\in X_1 \cup X_2\cup X_3} d_H(v)+ |H_3|
	\le |X_1\cup X_2|\cdot n + \sum_{v\in X_3} d_H(v)+ |H_3|
	\le 2p^{-2}n\log n + \sum_{v\in X_3} d_H(v).
\end{equation}
On the other hand, observe that for all vertices~$v\in X_3$ we have~$d_H(v) \le \eps pn$, due to the 
definition of~$X_2$.  Therefore, the number of missing edges is at least
\begin{equation}\label{eq:exploitMissingII}
\begin{split}
	\frac12 \sum_{v\in X_3} d_M(v)
	\ge & ~ \frac12  \sum_{v\in X_3} (d_H(v) + 5\eps pn)  \\
	\ge & ~ \frac12  \sum_{v\in X_3} (2d_H(v) + 4\eps pn)
	\ge \sum_{v\in X_3} d_H(v) + 2|X_3|  \eps pn.
\end{split}
\end{equation}
Now replace all horizontal edges from~$T$ with all missing edges in order to obtain a different 
triangle-free graph.  Comparing (\ref{eq:exploitMissingI}) and (\ref{eq:exploitMissingII}) yields that, 
to avoid a contradiction, $|X_3|$ must satisfy
\[
	2p^{-2}n\log n \ge 2|X_3|\eps pn 
	\;\Rightarrow\;
	|X_3| \le \eps^{-1}p^{-3}\log n \le p^{-4}\log n,
\]
whenever~$n$ is sufficiently large.

\medskip
\noindent {\bf (iv)} \quad 
Let~$v$ be a vertex in~$A_3$ with~$d_{H_3}(v) \ge p^{-6}\log n$ -- we handle vertices in~$B_3$ 
analogously.  Recall that~$A_3 = A\setminus (X_1 \cup X_2 \cup X_3)$. The definitions of the sets~$X_i$, 
together with the fact that~$|B| \ge (1-\eps)\frac n2$, imply that the number of neighbors of~$v$ in~$B$ 
is at least
\[
	d(T; v,B) = d(G; v,B) - d_M(v)
	\ge
	\frac{pn}3 - 6\eps pn
	\ge
	\frac{pn}4.
\]
Furthermore, note that the edges between the vertex sets~$\Gamma(T; v) \cap B$ and~$\Gamma(T; v) \cap A_3$ 
would form triangles in~$T$.  From (P3), we deduce that the number of missing edges in~$T$ is at least 
\begin{equation}\label{eq:setsZero}
	e(\Gnp; ~\Gamma(T; v) \cap B, ~\Gamma(T; v) \cap A_3)
	\ge
	\frac{p}2 \cdot \frac{pn}4 \cdot p^{-6}\log n 
	=
	\frac1{8}p^{-4}n\log n.
\end{equation}
Finally observe that the number of horizontal edges in~$T$ is, due to (i)--(iii), at most
\begin{equation}\label{eq:helperMissing}
	|H_3| + |X_1|\cdot n +  |X_2 \cup X_3|\cdot \frac54pn \le 2p^{-3}n\log n.
\end{equation}
But (\ref{eq:setsZero}) and (\ref{eq:helperMissing}) contradict the maximality of~$T$ -- we conclude that 
there are no vertices in~$A_3 \cup B_3$ with horizontal degree at least~$p^{-6}\log n$, i.e., (iv) is shown.

\medskip
\noindent {\bf (v)} \quad 
Let~$R$ be a matching of maximum cardinality in~$H_3$, and let~$m$ be the number of 
missing edges in~$T$. As in the previous arguments, the central idea is to derive a lower bound on~$m$ which 
contradicts the maximality of~$T$.

Let us assume without loss of generality that $e(T;\,A_3) \ge \frac12|H_3|$.  Using (iv), we readily obtain 
that~$|R| \ge \frac14|H_3|p^6(\log n)^{-1}$, as we can construct a matching by greedily removing edges 
from~$H_3$.  In order to estimate the number of missing edges, let~$e = \{u,v\}\in R$, where~$u,v \in A_3$. 
As~$u,v \not\in X_1$, the two vertices have at least~$\frac34p^2|B|\ge \frac{p^2n}4$ common neighbors in~$B$. 
Hence, for every~$e\in R$, there are at least~$\frac{p^2n}{4}$ missing edges.  As these edges are distinct 
for different edges in~$R$, we deduce that
\[
	m
	\ge
	|R| \cdot \frac{p^2n}{4}
	~\ge~ |H_3|\frac{p^6}{4\log n} \cdot \frac{p^2n}{4}
	~=~ |H_3| \cdot \frac{p^8}{16\log n}\;n.
\]
As in (\ref{eq:helperMissing}), we obtain that the number of horizontal edges is at most~$2p^{-3}n\log n$. 
In order to avoid a contradiction, $|H_3|$ must hence satisfy
\[
	|H_3| \cdot \frac{p^8}{16\log n}\;n \le 2p^{-3}n\log n,
\]
from which the claimed bound on~$|H_3|$ follows for sufficiently large~$n$.

\medskip
\noindent {\bf (vi)} \quad 
Let~$X = X_1 \cup X_2 \cup X_3$.  From our assumptions and (i)--(v), we know 
that~$|X| \le 2p^{-4}\log n$, if~$n$ is sufficiently large, and the number of horizontal edges 
in~$(A\cup B)\setminus X$ in~$T$ is less than~$\frac12p^{-12}\log^2n$.

To deduce the statement, we replicate the argument from (iv).  Let~$d$ be the maximum horizontal degree of a 
vertex in~$X$, and suppose~$d \ge 10 p^{-1}\log n$.  Furthermore, suppose that the maximum is attained at a 
vertex~$v\in A$.  Observe that~$v$ has at least~$d$ neighbors in~$B$, as otherwise~$\Pi$ would not have been 
optimal.  It follows that the number of missing edges is at least~$\frac{pd^2}{4}$, from (P3). 
On the other hand, the number of horizontal edges is at most~$|X|\cdot d  + \frac12p^{-12}\log^2n$. 
Therefore,~$d$ satisfies
\[
	\frac{pd^2}{4} \le 2p^{-4}\log n\cdot d + \frac12p^{-12}\log^2n
\]
which can only hold if~$d \le 3p^{-13/2}\log n$, and so the number of horizontal edges is at most
$p^{-12}\log^2 n$. This concludes the proof.
\end{proof}



\section{On Properties of (Near-)Optimal Bipartitions}\label{sec:cuts}

\renewcommand{\ne}{\overline{e}}
Before we proceed with the proof of Theorem \ref{thm:propertiesOfCuts}, we introduce two auxiliary tools. 
The first lemma is a statement about the number of non-edges between sufficiently large parts of the 
vertex set of the random graph~$\GnM$.  More precisely, for an ordered partition~$\Pi$ of the vertex set 
into two pairs of (sufficiently large) parts, we want to bound the probability that the number of 
\emph{non}-edges between the parts of the pairs is not near its expected value.  The result is not best 
possible, but it suffices for our purposes and keeps the calculations short.

\begin{lemma}\label{lem:numberNonEdgesCorrect}
Fix a natural number $s$ with~$n^{1/2} \le s \le \frac{n}{2}$, and let~$\Gamma = (V_1,W_1,V_2,W_2)$ be an 
ordered partition of~$[n]$ such that~$|V_1|+|V_2| = s$ and~$|W_1|,|W_2| \ge n/7$.  
Set $|\Gamma| := |V_1||W_1| + |V_2||W_2|$.  For a graph~$G$ with vertex set~$[n]$, let 
\begin{equation}\label{eq:nonEdgePartitions}
\ne(G;\, \Gamma) := |\Gamma| - e(G;\, V_1,W_1) - e(G;\, V_2,W_2).
\end{equation}
Take $M$ with~$n \le M \le \frac12\binom{n}{2}$ and set~$p := M/\binom{n}{2}$.  Then there is a 
constant~$C > 0$ such that
\[
	\Prob{\exists\Gamma: \big|\ne(\GnM; \Gamma) - (1-p)|\Gamma|\big| \ge C\cdot s^{-1/2}\cdot(1-p)|\Gamma|}
	\le e^{-n}.
\]
\end{lemma}

\begin{proof}
We show the analogous result for~$G_{n,p}$ and then use Pittel's inequality (\ref{eq:pittelsinequality}) 
to complete the proof.  For an ordered partition~$\Gamma$ with the above properties, define the event
\[
	\mathcal{E}_\Gamma ~:=~ \big|\ne(G_{n,p}; \Gamma) - (1-p)|\Gamma|\big| \ge C\cdot s^{-1/2}\cdot(1-p)|\Gamma|.
\]
Observe that, according to our assumptions, we have~$|\Gamma| \ge \frac{sn}{7}$.  A straightforward 
application of Lemma~\ref{len:cherny} then yields that we can choose~$C$ such 
that~$\Prob{\mathcal{E}_\Gamma} \le 20^{-n}$.  As there are at most~$4^n$ ways for choosing~$\Gamma$, the 
desired probability can be bounded above by $4^n \cdot 20^{-n} \ll n^{-1}e^{-n}$.
\end{proof}

The next proposition is a technical estimate for the probability that a trinomially distributed random 
variable of a special form deviates from its expectation.  The bound is tight up to the determined 
constant.  The proof is a straightforward application of Stirling's Formula, and we omit it.
\newcommand{\tp}{\alpha}
\newcommand{\tpN}{\tp N}

\begin{proposition}\label{prop:probabilityTrinomial}
Let~$\tp < \frac12$, and~$d \le \min\{\sqrt{\tpN}, \sqrt{(1-2\tp)N}\}$ such that~$2\tpN+d < N$. 
There is a constant~$C>0$ such that
\begin{equation}\label{eq:estimateMultinomial}
	\binom{N}{\tpN,\tpN+d} \tp^{2\tpN+d}(1-2\tp)^{(1-2\tp)N-d} \ge \frac{C}{\tpN}.
\end{equation}
\end{proposition}


With the above tools we are ready to prove the main result of this section.

\begin{proof}[Proof of Theorem \ref{thm:propertiesOfCuts}]
Recall that, for two bipartitions~$\Pi$ and~$\Pi'$ of the vertex set of a graph~$G$, we denote 
by~$\dist(\Pi, \Pi')$ the number of vertices in which~$\Pi$ and~$\Pi'$ differ, and we say that~$\Pi$ 
has \emph{gap}~$g$ (and write~$\gap(G;\, \Pi) = g$), if the number of edges across~$\Pi$ is 
exactly~$g$ less than the number of edges across an optimal bipartition.  Note that, for every pair 
of partitions~$\Pi$ and~$\Pi'$, we have $\dist(\Pi, \Pi') \le \frac{n}{2}$.

The central idea in our proof is to think of the random graph $G_{n,M}$ as being generated 
by adding edges one at a time, and to investigate how the random variable~$b(G_{n,M})$ changes as the 
number $M$ of edges increases.  More precisely, let~$t = t(n) > 0$ and consider the 
\emph{expected change}
\begin{equation}\label{eq:mlcEChange}
	\E{b(\GnMt) - b(\GnM)}.
\end{equation}
In order to obtain bounds for the above expression, we look at it from two different points of view: either 
removing~$t$ edges uniformly at random from~$G_{n, M+t}$ or adding~$t$ edges uniformly at random to~$\GnM$.

Suppose that we generate our random graph $G_{n,M}$ by first generating a uniform graph $G_{n,M+t}$, 
and then deleting~$t$ edges uniformly at random: we are interested in the expected change in $b$ on
the deletion of the $t$ edges.  Let~$\Pi^* = \Pi^*(\GnMt)$ be the canonical optimal bipartition of the 
graph~$\GnMt$.  Observe that, on deleting $t$ edges, the size of the optimal bipartition of the graph decreases 
by \emph{at most} the number of edges among the~$t$ deleted that go across~$\Pi^*$.  Hence, the expected 
decrease of the size of an optimal bipartition can be bounded from above by the expected proportion of edges 
of~$\GnMt$ across~$\Pi^*$.  Now we apply Proposition~\ref{prop:mlcGnM}, which yields that~$\GnMt$ has, with 
probability at least~$1-e^{-n}$, the property that all its optimal bipartitions have size at 
most~$\frac12(M+t) + \sqrt{4n(M+t)}$; we obtain 
\begin{equation}\label{eq:upper}
\begin{split}
	\E{b(\GnMt) - b(\GnM)}
	\le
	t \cdot \left( \frac{\frac12(M+t) + \sqrt{4n(M+t)}}{M + t} + e^{-n} \right)
	\le
	t\cdot\left(\frac12 + \sqrt{\frac{5n}{M}}\right).
\end{split}
\end{equation}

Recall that~$\dist(G;\Pi)$ denotes the distance of a bipartition $\Pi$ from the canonical optimal
bipartition~$\Pi^*(G)$.    
Let~$i^*:=i^*(n)$ be the minimum integer such that~$2^{i^*}s_0 \ge \frac{n}{4}$.  
For~$0\le i\le i^*$, define the event
\[
	\mathcal{P}_i ~:=~ \exists \Pi : \gap(G_{n,M}; \; \Pi) = r-1 \;	\text{ and }\;
		2^is_0\le \dist(G_{n,M};\Pi) \le 2^{i+1}s_0,
\]
and set~$p_i := \Prob{\mathcal{P}_i}$.  Note that the union of the events $\mathcal{P}_i$ is equivalent to 
the event
\[
	\mathcal{P} ~:=~ \exists \Pi : \gap(G_{n,M}; \; \Pi) = r-1 \; \text{ and }\;
		\dist(G_{n,M};\Pi) \ge s_0.
\]

We will show that, for sufficiently large~$n$, we have~$p_i \le \frac1{2^{i+1}\omega}$ for 
all~$0\le i \le i^*$.  The statement of the theorem then follows immediately:
\[
	\Prob{\mathcal{P}} 
	~\le~
	\sum_{i = 0}^{i^*(n)}p_i
	~\le~
	\frac1{2\omega} \sum_{i\ge 0} 2^{-i}
	~\le~ \frac1{\omega}.
\]

 From now on, we fix $i$ with~$0\le i \le i^*$, and set~$s_i := 2^is_0$ 
and~$t_i := r^2\frac{n(n-1)}{s_i(n-s_i)}$.  If we add~$t_i$ edges uniformly at random to~$G_{n,M}$, 
one of the following two events can occur:

\begin{itemize}
	\item[(i)] ~$\Pi^* = \Pi^*(G_{n,M}) = (A^*, B^*)$ remains one of the optimal bipartitions;
	\item[(ii)] a bipartition  different from~$\Pi^*$ ``overtakes''~$\Pi^*$, i.e., its size is 
		larger than the size of~$\Pi^*$ after having added~$t_i$ edges to~$G_{n,M}$.
\end{itemize}
Denote by~$\mathcal{I}_{t_i}$ the increase of the size of~$\Pi^*$, and by~$\io$ the indicator variable 
of the event that some bipartition different from~$\Pi^*$ overtakes~$\Pi^*$.  Due to linearity of 
expectation, we obtain
\begin{equation}\label{eq:expectedChangeLowerBound}
	\E{b(\GnMt_i) - b(\GnM)} \ge \E{\mathcal{I}_{t_i}} + \E{\io}.
\end{equation}
We first bound~$\E{\mathcal{I}_{t_i}}$.  For this, observe that~$\Pi^*$ increases by the number of edges 
among the~$t_i$ added that go across~$\Pi^*$; therefore, the expected increase of~$e(\GnM;\; \Pi^*)$
is~$t_i$ times the proportion of \emph{non}-edges across~$\Pi^*$ in~$G_{n,M}$.  Let~$C'>0$ be the constant 
that is guaranteed by Corollary \ref{cor:bmlcGnM}.  We apply Corollary \ref{cor:bmlcGnM} to obtain 
bounds for the minimum number of non-edges across any maximum bipartition of~$\GnM$, which hold with 
exponentially high probability; furthermore, due to our assumption~$M\le \frac12\binom{n}{2}$, we obtain 
that~$\binom{n}{2} - M \ge \frac14n^2$.  These facts imply that, for sufficiently large~$n$, we have
\newcommand{\Mb}{\widetilde{M}}
\begin{equation}\label{eq:lowerI}
	\E{\mathcal{I}_{t_i}}
	~\ge~
		(1 - 2e^{-n})
		\cdot
		t_i \cdot
		\frac{\frac12\left(\binom{n}{2} - M\right) - \sqrt{\frac{C'n^5}{M}}}{\binom{n}{2} - M}
	~\ge~
	t_i\cdot\left(\frac12 - \sqrt{\frac{20C'n}{M}}\right).
\end{equation}
\newcommand{\EBAL}{\mathcal{E}_{\mathcal{B}al}}%
\newcommand{\tmA}{\overline{e}_\mathcal{X}}%
\newcommand{\tmB}{\overline{e}_\mathcal{Y}}%
\newcommand{\mA}{\mathcal{X}}%
\newcommand{\mB}{\mathcal{Y}}%
\newcommand{\EXY}{\mathcal{E}_{\mA\mB}}%
\newcommand{\ED}{\mathcal{E}_{edges}}%
\newcommand{\EL}{\mathcal{E}_{\Lambda}}%

Next we bound~$\E{\io}$ from below; our first step is to specify two properties that we want of
the graph $G_{n,M}$.  Denote by~$\EBAL$ the event that there is no bipartition $\Pi =(A,B)$
with $\gap(\GnM;\; \Pi) \le r-1$, and $\min (|A|, |B|) < \frac{n}{2} - 4n^{\frac34}p^{-\frac14}$.  
Observe that Proposition \ref{prop:mlcIsBalanced} implies~$\Prob{\GnM \in \EBAL}\ge 1 - e^{-n}$, with 
room to spare.

Let~$\tilde{C}$ be the constant guaranteed by Lemma~\ref{lem:numberNonEdgesCorrect} and let $\ED$ 
be the event that, for all $s$ with $n^{1/2} \le s \le \frac n2$, and all ordered partitions
$\Gamma=(V_1,W_1,V_2,W_2)$ with $|V_1| + |V_2| = s$, $|W_1|, |W_2| \ge n/7$, we have
$\big|\ne(G_{n,M};\, \Gamma) - (1-p)|\Gamma|\big| \le \tilde{C}s^{-1/2}(1-p)|\Gamma|$.  
By Lemma~\ref{lem:numberNonEdgesCorrect}, we have $\Prob{\GnM \in \ED} \ge 1-ne^{-n}$.  

Suppose that $\EBAL$ and $\ED$ occur.  Suppose also that ${\mathcal P}_i$ occurs, and fix 
some bipartition~$\Pi = (A,B)$ of~$G_{n,M}$ with the properties~$\gap(\GnM;\;\Pi)=r-1$ 
and~$\dist(\Pi,\Pi^*) =: s$ with~$s_i \le s\le 2s_i$.  
Note that, from $\EBAL$, we do then have 
\begin{equation}\label{eq:PiPisbalanced}
        |A^*|, |B^*|, |A|, |B| \ge \frac{n}{2} - 4n^\frac34 p^{-\frac14}.
\end{equation}
In order to avoid ambiguities, we label the sets so that~$|A^*| \ge |B^*|$ and
\[
	s = |A \cap B^*| + |A^* \cap B| \le |A \cap A^*| + |B \cap B^*|.
\]
Note that if we add~$t_i$ edges uniformly at random to~$\GnM$, edges added between~$A^* \cap B$ 
and~$A^* \cap A$, and between~$B^* \cap A$ and~$B^* \cap B$, contribute to~$e(\GnM; \, \Pi)$, while 
edges added between~$A^* \cap B$ and~$B^* \cap B$, and between~$B^* \cap A$ and~$A^* \cap A$, 
contribute to~$e(\GnM; \, \Pi^*)$.  All other added edges contribute simultaneously to 
both~$e(\GnM;\, \Pi^*)$ and~$e(\GnM;\, \Pi)$, or to neither of them.  This motivates the definition of 
the following two derived ordered partitions
\[
	\Gamma_\mA := (A^* \cap B, A^* \cap A, B^* \cap A, B^* \cap B)
	~\text{ and }~
	\Gamma_\mB := (A^* \cap B, B^* \cap B, B^* \cap A, A^* \cap A).
\]
Additionally, let~$|\Gamma_\mA|$,~$\ne(\GnM; \Gamma_\mA)$, and similarly~$|\Gamma_\mB|$ 
and~$\ne(\GnM; \Gamma_\mB)$, be defined as in Lemma~\ref{lem:numberNonEdgesCorrect}; 
see~(\ref{eq:nonEdgePartitions}).  As discussed above, the motivation for this definition is that the 
event~$\mathcal{O}$ will occur if, out of the~$t_i$ edges added to~$G_{n,M}$, the number of edges added 
among those counted by~$\ne(\GnM; \Gamma_\mA)$ is at least~$r$ plus the number of edges added among those 
counted by~$\ne(\GnM; \Gamma_\mB)$.  We shall denote this event by~$\EXY$: note that, once we have fixed 
$G_{n,M}$ and $\Pi$, this event depends solely on the choice of the $t_i$ added edges.  


Observe that, due to our assumptions, we have
\begin{equation}\label{eq:atLeastSqrt}
	|A^* \cap B| + |B^* \cap A| = s \ge s_0 \ge r^4\sqrt{np^{-1}} \ge \sqrt{n}.
\end{equation}
Furthermore, due to (\ref{eq:PiPisbalanced}), we have for sufficiently large~$n$ that
\[
	\big| |A^* \cap B| - |B^* \cap A|\big| = \big| |A^*| + |B| - (|A^*\cup B| + |A\cap B^*|) \big|
	= \big| |A^*| + |B| - n \big| \le 8n^{3/4}p^{-1/4},
\]
which implies
\begin{equation}\label{eq:deltaParts}
	|A^* \cap B|, |A \cap B^*|
	\le
	\frac{s}{2} + 4n^{3/4}p^{-1/4}
	\stackrel{(s \le n/2)}{\le}
	\frac{n}{4} + 4n^{3/4}p^{-1/4}
	\le
	\frac{n}{3}.
\end{equation}
Hence, whenever~$n$ is sufficiently large, we obtain that
\[
	|A\cap A^*| = |A| - |A\cap B^*| \ge \frac{n}{2} - 4n^{3/4}p^{-1/4} - \frac{n}{3} \ge \frac{n}{7},
\]
and similarly,~$|B\cap B^*| \ge \frac{n}7$.  Together with (\ref{eq:atLeastSqrt}), and the assumption
that the event $\ED$ does not occur, this implies that, for $\mathcal{S} \in \{ \mA, \mB\}$, 
\[
\big|\ne(G;\, \Gamma_\mathcal{S}) - (1-p)|\Gamma_\mathcal{S}|\big| \le
                               \tilde{C}s^{-1/2}(1-p)|\Gamma_\mathcal{S}|.
\]
The above discussion yields
\begin{equation}\label{eq:bringprsIntoGame}
	\E{\io} = \Prob{{\mathcal O}}
	\ge
	\Prob{\mathcal{P}_i \wedge \EBAL \wedge \ED}
		\cdot
		\condProb{\EXY}{\mathcal{P}_i \wedge \EBAL \wedge \ED}.
\end{equation}
As discussed above, we can estimate 
\[
	\Prob{\mathcal{P}_i \wedge \EBAL \wedge \ED}
	= p_i - \Prob{\mathcal{P}_i \setminus (\EBAL \wedge \ED)}
	\ge p_i - \Prob{\overline{\EBAL \wedge \ED}}
	\ge p_i - 2ne^{-n}.
\]
Our aim is to show that there is a constant~$c_1>0$ such that
\begin{equation}\label{eq:expectedChangeFinalBound}
	\condProb{\EXY}{\mathcal{P}_i \wedge \EBAL \wedge \ED}
	~\ge~
	c_1 \cdot r^{-2}.
\end{equation}
Before proving this, we explain how this will complete the proof of the theorem.  From (\ref{eq:upper}), 
where we set~$t = t_i$, as well as (\ref{eq:expectedChangeLowerBound}), (\ref{eq:lowerI}), 
(\ref{eq:bringprsIntoGame}) and the above estimate, we obtain
\[
	t_i\cdot\left(\frac12 - \sqrt{\frac{20C'n}{M}}\right)
	+ (p_i - 2ne^{-n})\cdot	\frac{c_1}{r^2}
	~\le~
	t_i\cdot\left(\frac12 + \sqrt{\frac{5n}{M}}\right),
\]
Using also that $s \le \frac{n}{2}$, and hence~$t_i = r^2\frac{n(n-1)}{s_i(n-s_i)} \le \frac{2r^2n}{s_i}$,
this implies that we can choose a constant~$c_2$ such that, for sufficiently large~$n$,
\[
	p_i
	~\le~ 
	c_2 \cdot t_i\cdot \sqrt{\frac{n}{M}} \cdot r^2
	~\le~
	4c_2 \cdot \frac{r^4n}{s_i} \sqrt{\frac{1}{pn}}
	~\le~
	\frac{8c_2 \cdot r^4}{2^{i+1}\cdot C\cdot\omega\cdot r^4\cdot\sqrt{np^{-1}}} \sqrt{\frac{n}{p}}
	~\le~
	\frac{1}{2^{i+1}\omega},
\]
where we set~$C := 8c_2$.

\newcommand{\EX}{\overline{E}_\mathcal{X}}
\newcommand{\EY}{\overline{E}_\mathcal{Y}}
Let~$q := \condProb{\EXY}{\mathcal{P}_i \wedge \EBAL \wedge \ED}$; to complete the proof we estimate~$q$ 
so as to obtain~\eqref{eq:expectedChangeFinalBound}.  Suppose then that the graph $G_{n,M}$ does
satisfy $\EBAL$, $\ED$ and ${\mathcal P}_i$, and that we have fixed a bipartition $\Pi$ as above.   

For brevity, let~$\tmA := \ne(\GnM;\, \Gamma_\mA)$ and~$\tmB := \ne(\GnM;\, \Gamma_\mB)$.  Furthermore, let~$\EX$ 
denote the set of edges counted in~$\tmA$, and define similarly~$\EY$.  In order to obtain a lower bound 
for~$q$, it suffices to consider the event that precisely~$x_i:=\frac{s(n-s)}{n(n-1)}t_i$ edges out of~$\EY$ 
are added to~$\GnM$ and~$x_i+r$ out of~$\EX$.  Suppose that we add the~$t_i$ edges one after the other -- then 
there are~$\binom{t_i}{x_i+r+1,x_i}$ ways to choose the points in time at which those edges are taken 
from~$\EX$ and~$\EY$.  Furthermore, the probability that an edge of~$\EX$ is added to~$\GnM$ is at 
least~$\frac{\tmA - t_i}{\binom{n}{2}-M}$, and similarly, the probability that an added edge is contained 
in~$\EY$ can be bounded from below by~$\frac{\tmB - t_i}{\binom{n}{2}-M}$.  Moreover, the probability that 
an added edge belongs neither to~$\EX$ nor~$\EY$ is at least~$1 - \frac{\tmA+\tmB}{\binom{n}{2}-M-t_i}$. 
Putting all this together yields
\begin{equation}\label{eq:EXY}
	q
	\ge
	\binom{t_i}{x_i,x_i+r}
	\cdot
	\left(\frac{\tmA-t_i}{\binom{n}{2} - M}\right)^{x_i+r}
	\left(\frac{\tmB-t_i}{\binom{n}{2} - M}\right)^{x_i}
	\left(1 - \frac{\tmA + \tmB}{\binom{n}{2} - M - t_i}\right)^{t_i-2x_i-r}.
\end{equation}
We now simplify the above expression so that we can apply Proposition~\ref{prop:probabilityTrinomial} to 
obtain~\eqref{eq:expectedChangeFinalBound}.  As a first step, we derive tight bounds for the 
quantities~$|\Gamma_\mathcal{X}|$ and~$|\Gamma_\mathcal{Y}|$.  
Let~$s^{(1)} := \card{A^* \cap B}$, $s^{(2)} := \card{B^* \cap A}$, and note that~$s = s^{(1)} + s^{(2)}$. 
Furthermore, define the quantity~$\delta := s^{(1)} - \frac{s}{2}$. Observe that
\[
	\card{\Gamma_\mB}
	=
	s^{(1)}(|B^*| - s^{(2)}) + s^{(2)}(|A^*| - s^{(1)})
	=
	\frac{s(n-s)}{2} - \underbrace{\delta(|A^*| - |B^*|  - 2\delta)}_{:=\Delta}.
\]
Due to the event~$\EBAL$, we 
have~$\big||A^*| - |B^*|\big| \le 8n^{3/4}p^{-1/4}$ and~$|s^{(1)} - s^{(2)}| \le 8n^{3/4}p^{-1/4}$.  This 
implies that
\begin{equation*}
	|\delta| \le \min\left\{\frac{s}{2}, 4n^{3/4}p^{-1/4}\right\}
	~\text{ and hence }~
	|\Delta| \le 32\min\left\{sn^{3/4}p^{-1/4}, n^{3/2}p^{-1/2}\right\}.
\end{equation*}
 From~$s \ge s_0 = C\omega r^4\sqrt{np^{-1}}$, we can deduce
\begin{equation}\label{eq:estimateGammaX}
	|\Gamma_\mB|
	= \frac{s(n-s)}{2} \Big(1 + \Theta\big(\min\big\{(pn)^{-1/4}, {\omega^{-1} r^{-4}}\big\}\big)\Big)
	\stackrel{(r\ll (pn)^{1/8})}{=} \frac{s(n-s)}{2}\big(1 + o(r^{-2})\big).
\end{equation}
Now note that~$|\Gamma_\mA| + |\Gamma_\mB| = s(n-s)$ implies~$\card{\Gamma_\mA} = \frac{s(n-s)}{2} + \Delta$; hence 
the above statement is also valid for~$|\Gamma_\mA|$.

Before we proceed let us make some auxiliary calculations, which we will need in the remainder.  First,
\begin{equation*}
	x_i
	= \frac{s(n-s)}{n(n-1)} \cdot t_i 
	= \frac{s(n-s)}{n(n-1)} \cdot r^2\frac{n(n-1)}{s_i(n-s_i)}
	\stackrel{(s_i\le s \le 2s_i)}{\le} 2r^2.
\end{equation*}
Due to (\ref{eq:estimateGammaX}) we may assume, for sufficiently large~$n$, 
that~$|\Gamma_\mA|, |\Gamma_\mB| \ge \frac{sn}{5}$, which implies, together with the event~$\ED$, that we 
have~$\tmA \ge (1 - \tilde{C}s^{-1/2})(1-p)|\Gamma_\mA| \ge \frac{sn}{11}$.  With~$s_i\le s\le\frac{n}{2}$ 
and~$s_i \ge s_0 \ge \omega r^4 \sqrt{np^{-1}}$ we can deduce, if~$n$ is large enough, that
\begin{equation}\label{eq:helper}
	\frac{t_ix_i}{\tmB}
	\le \frac{11}{sn}
		\cdot \frac{r^2n(n-1)}{s_i(n-s_i)}
		\cdot 2r^2
	\le \frac{22r^4}{s_i^2}
	= o(1)
	~~\text{ and }~~
	s^{-\frac12}x_i
	\le \frac{2r^2}{(\omega r^4\sqrt{np^{-1}})^{\frac12}}
	\le n^{-\frac14}.
\end{equation}
Now we proceed with simplifying (\ref{eq:EXY}).  Using the inequality~$1-y\ge e^{-2y}$, which is valid for 
e.g.~$0\le y \le \frac12$, the third term on the right-hand side of (\ref{eq:EXY}) can be estimated for 
sufficiently large~$n$ by
\begin{equation}\label{eq:estimateX}
\begin{split}
	\left(\frac{\tmB-t_i}{\binom{n}{2} - M}\right)^{x_i}
	&\ge \left(\frac{\tmB}{\binom{n}{2} - M}\right)^{x_i} \cdot \left(1 - \frac{t_i}{\tmB}\right)^{x_i} \\
	&\stackrel{(\ED)}{\ge}
		\left(\frac{(1 - \tilde{C}s^{-1/2})(1-p)|\Gamma_\mB|}{(1-p)\binom{n}{2}}\right)^{x_i}
		\cdot e^{-\frac{2t_ix_i}{\tmB}} \\
	&\,\stackrel{(\ref{eq:helper})}{\ge}
		\left(\frac{|\Gamma_\mB|}{\binom{n}{2}}\right)^{x_i}
		\cdot
		\left(1 - \tilde{C}s^{-1/2}\right)^{x_i}
		\cdot \frac12 \\
	&\text{\hspace{-0.4cm}}\stackrel{(\ref{eq:estimateGammaX}), (\ref{eq:helper})}{\ge}
		\left(\frac{s(n-s)}{2\binom{n}{2}}\right)^{x_i}
		(1 + o(r^{-2}))^{x_i}
		\cdot e^{-\frac{4\tilde{C}}{n^{1/4}}} \cdot \frac12\\
	&\ge \frac14\left(\frac{s(n-s)}{2\binom{n}{2}}\right)^{x_i}.
\end{split}
\end{equation}
Precisely the same calculation yields that the second term of the right-hand side of~\eqref{eq:EXY} is, 
for sufficiently large~$n$, at least
\begin{equation}\label{eq:estimateY}
	\left(\frac{\tmA-t_i}{\binom{n}{2} - M}\right)^{x_i+r}
	\ge \frac14 \left(\frac{s(n-s)}{2\binom{n}{2}}\right)^{x_i+r}.
\end{equation}
In order to simplify the last term in (\ref{eq:EXY}), first recall that~$|\Gamma_\mA| + |\Gamma_\mB| = s(n-s)$. 
Due to the event~$\ED$, we have~$\tmA + \tmB \le (1 + \tilde{C}s^{-1/2})(1-p)s(n-s)$. 
Abbreviate~$m := \binom{n}{2} - M = (1-p)\binom{n}{2}$, and observe that we have~$m\ge \frac15n^2$.  
This yields
\[
\begin{split}
	\left(1 - \frac{\tmA + \tmB}{m - t_i}\right)^{t_i-2x_i-r}
	&\ge
	\left(1 - \frac{(1 + \tilde{C}s^{-1/2})(1-p)s(n-s)}
	{(1 - \frac{5t_i}{n^2})(1-p)\binom{n}{2}}\right)^{t_i-2x_i-r}.
\end{split}
\]
Note that, due to~$r \ll (pn)^{1/8} \le n^{1/8}$, we have~$\frac{5t_i}{n^2} \ll \frac{1}{n}$, which 
implies~${s^{-1/2}}(1 - \frac{5t_i}{n^2})^{-1} \le 2s^{-1/2}$, whenever~$n$ is sufficiently large.  
Thus, using that the estimate~$1-y \ge e^{-2y}$, valid for sufficiently small $y$, we obtain for large~$n$ 
that
\begin{equation}\label{eq:estimateXY}
\begin{split}
	\left(1 - \frac{\tmA + \tmB}{m - t_i}\right)^{t_i-2x_i-r}
	&\ge
		\left(1 - \frac{s(n-s)}{\binom{n}{2}}\right)^{t_i-2x_i-r}
		\left(1 - 5\tilde{C}s^{-1/2}\frac{s(n-s)}{\binom{n}{2}}\right)^{t_i} \\
	&\ge
		\left(1 - \frac{s(n-s)}{\binom{n}{2}}\right)^{t_i-2x_i-r}
		e^{-20\tilde{C}s^{-1/2}x_i} \\
	&\hspace{-2mm}\stackrel{\eqref{eq:helper}}{\ge}
		\frac12\left(1 - \frac{s(n-s)}{\binom{n}{2}}\right)^{t_i-2x_i-r}.
\end{split}
\end{equation}
By combining (\ref{eq:estimateX}), (\ref{eq:estimateY}), and (\ref{eq:estimateXY}), we obtain from 
(\ref{eq:EXY}) that \[
	q ~\ge~
	\frac1{32}\binom{t_i}{x_i,x_i+r}
	\left(\frac{s(n-s)}{2\binom{n}{2}}\right)^{2x_i+r}
	\left(1 - \frac{s(n-s)}{\binom{n}{2}}\right)^{t_i-2x_i-r}.
\]
Now it can be easily checked that we can apply Proposition \ref{prop:probabilityTrinomial} to estimate 
the above expression, where we set~$N := t_i$,~$\alpha := \frac{s(n-s)}{n(n-1)}$,~$d := r$. 
Indeed,~$\alpha$ reaches its maximum value when~$s = \frac{n}{2}$, and we obtain for~$n\ge 4$ 
that~$\alpha \le \frac13 < \frac12$; furthermore,~$t_i = r^2\frac{n(n-1)}{s_i(n-s_i)}$ 
implies~$d \le \sqrt{\min\{\alpha, 1-2\alpha\}N}$.  Hence, if we denote by~$C''$ the constant defined by 
Proposition~\ref{prop:probabilityTrinomial}, and use that~$s_i \le s\le 2s_i$, we obtain
\[
	q
	~\ge~	\frac1{64}\cdot\frac{C''}{x_i}
	~\ge~ \frac{C''}{64} \cdot \frac{1}{r^2}\frac{s_i(n-s_i)}{s(n-s)}
	~\ge~ c_1r^{-2},
\]
which is precisely (\ref{eq:expectedChangeFinalBound}), if we choose~$c_1$ appropriately.  
This completes the proof.
\end{proof}


\section{Proof of Theorem \ref{thm:extremalMain}}\label{sec:FKGproof}

\newcommand{\BalPi}{\mathcal{B}al}

\newcommand{\Fnpl}{\mathcal{F}_\ell(G_{n,p})}
\newcommand{\Gnp}{G_{n,p}}
\newcommand{\Eall}{\mathcal{E}}
\newcommand{\Eone}{\mathcal{E}_1}
\newcommand{\Etwo}{\mathcal{E}_2}
\newcommand{\Bone}{\mathcal{B}_1}
\newcommand{\Btwo}{\mathcal{B}_2}
\newcommand{\Ball}{\mathcal{B}}
\newcommand{\notB}{\overline{\mathcal{B}}}

Let~$n^{-1/250} \le p \le \frac12$. Furthermore, define the functions
\[
	r_0 = r_0(p,n) := p^{-12}\log^2n ~~\text{ and }~~ s(r) := n^{2/3} \cdot (r+1)^4,
\]
and abbreviate~$s_0 := s(2r_0)$. Before we continue, we need to introduce some notation. 
Let~$G\in\mathcal{G}_n$ be a graph with vertex set~$[n] := \{1,\dots,n\}$. A bipartition of~$[n]$ is said 
to be \emph{balanced} if it is contained in the set
\[
	\BalPi_{n} := \left\{
			(A,B) ~|~
			\card{|A| - \frac{n}{2}} \le \frac{n}{100}
			~\text{ and }~
			\card{|B| - \frac{n}{2}} \le \frac{n}{100}
		\right\}.
\]
Recall that~$\mathcal{T}(G)$ denotes the set of maximum-cardinality triangle-free subgraphs of a graph~$G$. 
Define the two ``bad'' events
\begin{equation}\label{eq:badEvents}
\begin{split}
	\Bone := \Big\{
		G \in \mathcal{G}_n ~|~ &\text{$\exists T\in\mathcal{T}(G) :$ 
		there is no balanced bipartition~$\Pi$ } \\
		&~~~~~~~\text{with at most~$r_0$ edges of~$T$ inside~$\Pi$}
		\Big\}, \\
	\Btwo := \Big\{
		G \in \mathcal{G}_n ~|~ &\exists\Pi :
		\gap(G; \;\Pi) \le 2r_0
		\text{ and }
		\dist(G; \; \Pi) \ge \frac{s\big(\gap(G; \;\Pi)\big)}{2}
		\Big\}.
\end{split}
\end{equation}
Moreover, let~$\Ball = \Bone \cup \Btwo$ and define the ``good'' 
event~$\notB = \mathcal{G}_n \setminus \Ball$. 

In the sequel, we estimate the probability that there is a~$T\in\mathcal{T}(\Gnp)$ that is not bipartite. 
Observe that this implies the existence of a bipartition~$\Pi_T = \Pi = (A,B)$ of the vertex set with the 
property that we can obtain~$T$ from the subgraph~$E(\Gnp;\, \Pi)$ of the random graph by 
removing~$t > 0$ edges, and adding some number at least~\mbox{$t + \gap(\Gnp; \, \Pi)$} edges 
of~$E(\Gnp)\setminus E(\Gnp;\, \Pi)$.  Accordingly, for a fixed~$\Pi$ and any set~$S$ of edges 
inside~$\Pi$, let
\begin{equation}\label{eq:eventE}
\begin{split}
	\Eall(\Pi, S) := \big\{
			G \in \mathcal{G}_n ~|~
			\exists X\subseteq E(G;\, \Pi) : 
			(E(G;\, \Pi)\setminus X) & \cup S \text{ is triangle-free } \\
			&\text{and }  |S| - |X| \ge \gap(G;\, \Pi)
			\big\}.
\end{split}
\end{equation}
Now let us assume that~$\Gnp\in\notB$.  If there is a~$T\in\mathcal{T}(\Gnp)$ which is not bipartite, then 
there exists a balanced partition~$\Pi = (A,B)$, which can be ``enhanced'' by at least one and at 
most~$r_0$ edges, with both endpoints in~$A$ or~$B$, and possibly by removing at most~$r_0$ edges with one 
endpoint in~$A$ and one in~$B$, such that we obtain a triangle-free graph.  The above definition of the 
event~$\Eall(\Pi, S)$ thus implies that we have
\begin{equation}\label{eq:getRidOfBadEvents}
\begin{split}
	\Pr[\exists T\in\mathcal{T}&(\Gnp): T \text{ is not bipartite}] \\
	&\le
		\Prob{\Gnp\in\Ball}
		+ \sum_{\substack{\Pi=(A,B)\in\BalPi_n \\ S\in \binom{A}{2}\cup\binom{B}{2}, 
		1 \le |S|\le r_0}}\Prob{\Gnp\in(\Eall(\Pi,S) \cap \notB)}.
\end{split}
\end{equation}

Let~$C$ be the constant guaranteed by Theorem~\ref{thm:propertiesOfCuts}.  By applying 
Lemma~\ref{lem:almostlpartite} and Theorem~\ref{thm:propertiesOfCuts}, for all~$r\in[1,2r_0+1]$ 
and~$\omega := \frac{1}{2C}n^{1/6-1/250}$, the first term on the right hand side 
of~(\ref{eq:getRidOfBadEvents}) can be bounded as follows for sufficiently large~$n$:
\begin{equation}\label{eq:boundBadEvents}
\begin{split}
	\Prob{\Gnp\in\Ball}
		&\;\le\;
		\Prob{\Gnp\in\Bone} + \Prob{\Gnp\in\Btwo} \\
		&\;\le\; o(1) + \sum_{r=1}^{2r_0+1} \frac{1}{\omega}
		\;=\; o(1) + \frac{2n^{12/250}\log^2n + 1}{\omega}
		\;=\; o(1).
\end{split}
\end{equation}
In the remainder of the proof, we will bound the sum on the right hand side of (\ref{eq:getRidOfBadEvents}). 
In particular, we show that the probability of the joint event~$(\Gnp\in\Eall(\Pi,S)) \wedge (\Gnp\in\notB)$ 
can be estimated by the probability that two appropriately defined events~$\Eone$ and~$\Etwo$ occur.  This 
will allow us to use the FKG inequality in order to get a sufficient upper bound for the sum 
in~\eqref{eq:getRidOfBadEvents}.  For a fixed~$\Pi = (A,B)\in\BalPi_n$ and a non-empty 
set~$S \subseteq\binom{A}{2}\cup\binom{B}{2}$ of edges, let
\begin{equation}\label{eq:E1E2}
\begin{split}
	\Eone(\Pi) &= \Big\{
			G \in \mathcal{G}_n ~|~
			\gap(G;\;\Pi) \le r_0
			\text{ and }
			\forall\Pi': \big(\gap(G;\; \Pi,\Pi') \le r_0 \Rightarrow \dist(\Pi,\Pi') \le s_0\big)
			\Big\}, \\
	\Etwo(\Pi,S) &= \Big\{
			G \in \mathcal{G}_n ~|~
			\exists X\subset E(G;\Pi) :
			(E(G; \Pi)\setminus X)\cup S
			\text{ is~triangle-free, and }
			|X| \le |S|
			\Big\}.
\end{split}
\end{equation}
The definition of~$\Eone$ may seem overly complicated, and one could think that an event of the type 
``there is no partition~$\Pi$ with~$\gap(G;\,\Pi) \le 2r_0$ and~$\dist(G;\,\Pi) \ge s_0$'', which would 
follow directly from the definition (\ref{eq:badEvents}) of~$\Btwo$, could be sufficient as well. 
It turns out (see Proposition~\ref{prop:FKGApplication}) that this is in fact a delicate point of the proof: 
we \emph{need} to relax the event, so that it becomes an increasing function in an appropriately defined 
partial ordering of all graphs with~$M$ edges.  As a consequence, we may use the FKG inequality to get a 
sufficient estimate.

First we show that in fact we can bound the probability from (\ref{eq:getRidOfBadEvents}) with the joint 
probability of the events~$\Eone$ and~$\Etwo$.

\begin{proposition}\label{prop:tranferToE1E2}
For all~$\Pi = (A,B)\in\BalPi_n$ and~$S\subseteq\binom{A}{2} \cup \binom{B}{2}$ with~$1\le |S|\le r_0$, 
\begin{equation}\label{eq:transferToE1E2}
	\Prob{(\Gnp\in\Eall(\Pi,S)) \wedge (\Gnp\in\notB)}
	\le
	\Prob{(\Gnp\in\Eone(\Pi)) \wedge (\Gnp\in\Etwo(\Pi,S))}.
\end{equation}
\end{proposition}

\begin{proof}
Suppose that the event~$(\Gnp\in\Eall(\Pi,S)) \wedge (\Gnp\in\notB)$ occurs.  From (\ref{eq:eventE}) we 
deduce that~$\gap(\Gnp; \; \Pi) \le |S| \le r_0$; hence,~$\Gnp\not\in\Btwo$ implies 
that~$\dist(\Gnp; \; \Pi) \le \frac{s(r_0)}{2}$.  Now consider any partition~$\Pi'$ different from~$\Pi$, 
which fulfills~$\gap(\Gnp; \; \Pi,\Pi') \le r_0$.  Clearly, we have that~$\gap(\Gnp;\;\Pi') \le 2r_0$, and 
consequently ~$\Gnp\not\in\Btwo$ implies~$\dist(\Gnp; \; \Pi') \le \frac{s(2r_0)}{2}$.  
As the distance of~$\Pi$ and~$\Pi'$ is at most the sum of their distances from the canonical maximum bipartition, 
we obtain
\[
	\dist(\Pi, \Pi') \le \dist(\Gnp;\, \Pi) + \dist(\Gnp;\, \Pi') \le s(2r_0) = s_0.
\]
That is, we have~$\Gnp \in \Eone(\Pi)$.  The event~$\Etwo(\Pi,S)$ is easily seen to hold simultaneously, 
as it is a relaxation of the condition in (\ref{eq:eventE}).
\end{proof}

The next proposition states that we can bound the probability from (\ref{eq:transferToE1E2}) by the 
product of the probabilities of the events~$\Eone$ and~$\Etwo$.  Its proof consists of the definition of an 
appropriate distributive lattice on graphs with respect to bipartitions of the vertex set, and a subsequent 
application of the FKG inequality.

\begin{proposition}\label{prop:FKGApplication}
For all ~$\Pi = (A,B)\in\BalPi_n$ and~$S\subset\binom{A}{2}\cup\binom{B}{2}$ with~$1\le |S|\le r_0$
\[
	\Prob{(\Gnp\in\Eone(\Pi)) \wedge (\Gnp\in\Etwo(\Pi,S))}
	\le
	\Prob{\Gnp\in\Eone(\Pi)} \cdot \Prob{\Gnp\in\Etwo(\Pi,S)}.
\]
\end{proposition}

\begin{proof}
\newcommand{\mL}{\mathcal{L}}
For the proof we need a variant of the FKG inequality which we shall now state; a far more general treatment 
of the topic can be found in \cite{book:as00}.  A \emph{lattice} is a partially ordered set~$(S, \le)$ (with 
ground set~$S$ and a partial order~$\le$ on~$S$) in which every two elements~$x$ and~$y$ have a unique 
minimal upper bound and a unique maximal lower bound, which we denote by~$x\vee y$ and~$x\wedge y$ 
respectively.  The lattice~$\mL$ is called \emph{distributive} if, for all~$x,y,z\in\mL$, we have
\begin{equation}\label{eq:distrLattice}
	x \wedge (y \vee z) = (x\wedge y) \vee (x\wedge z).
\end{equation}
A function~$f:S\rightarrow \mathbb{R}^+$ is called \emph{log-supermodular} if, for all~$x,y\in S$, 
\begin{equation}\label{eq:submodular}
	f(x) f(y) \le f(x\vee y)f(x\wedge y).
\end{equation}
A function~$f:S \rightarrow \mathbb{R}^+$ is called \emph{increasing} if~$x\le y$ implies~$f(x) \le f(y)$, 
and \emph{decreasing} if~$x\le y$ implies~$f(x) \ge f(y)$.  With these definitions, we can state a 
probabilistic version of the well-known FKG inequality: if~$\mL = (\Omega, \le)$ is a finite distributive 
lattice,~$f:\Omega\rightarrow\mathbb{R}^+$ is an increasing function,~$g:\Omega\rightarrow\mathbb{R}^+$ is a 
decreasing function, and~$\mu:\Omega\rightarrow\mathbb{R}^+$ is a log-supermodular probability measure 
on~$\Omega$, then we have
\begin{equation}\label{eq:myfkg}
	\E{f \cdot g} \le \E{f} \cdot \E{g}.
\end{equation}

We now prove Proposition \ref{prop:FKGApplication}.  In order to do so, we fix some~$\Pi$ and~$S$.  In order 
to apply inequality (\ref{eq:myfkg}), we define the following partial ordering on graphs with vertex set $[n]$.  
For two graphs~$G$ and~$H$ let
\begin{equation}\label{eq:ordering}
	G \le_\Pi H ~:\Leftrightarrow~ E(G; \Pi) \subseteq E(H; \Pi)							 ~\text{ and }~
	(E(G)\setminus E(G; \Pi)) \supseteq (E(H)\setminus E(H; \Pi)).
\end{equation}
Intuitively, a graph~$G$ is ``smaller'' than a graph~$H$ with respect to~$\le_\Pi$, if it has fewer edges 
across~$\Pi$ \emph{and} simultaneously more edges inside~$\Pi$. One easily checks that
for any pair of graphs~$G$ and~$H$ the unique minimal upper bound of~$G$ and~$H$ is given by 
\[
 G\vee H = \big( E(G; A) \cap E(H; A) \big)
       ~\cup~
      \big( E(G; B) \cap E(H; B) \big)
       ~\cup~ 
      \big( E(G; \Pi) \cup E(H; \Pi)\big),
\]
while the unique maximal lower bound is given by
\[
 G\wedge H := \big( E(G; A) \cup E(H; A) \big)
        ~\cup~
       \big( E(G; B) \cup E(H; B) \big)
        ~\cup~ 
       \big( E(G; \Pi) \cap E(H; \Pi)\big).
\]
It can easily be verified that these operators are distributive, i.e., they fulfill (\ref{eq:distrLattice}).
For our probability space, we use that of~$\Gnp$, i.e., for any~$G\in\mathcal{G}_n$, we 
set~$\mu(G) := \Prob{\Gnp = G} = p^{e(G)}(1-p)^{\binom{n}{2} - e(G)}$.  An easy calculation yields that 
\[
\begin{split}
	\mu(G)\mu(H)
	&= p^{e(G) + e(H)}(1-p)^{2\binom{n}{2} - e(G) - e(H)} \\
	&= p^{e(G\vee H) + e(G\wedge H)}(1-p)^{2\binom{n}{2} - e(G\vee H) - e(G\wedge H)}
	= \mu(G \vee H)\mu(G\wedge H),
\end{split}
\]
i.e.,~$\mu$ is log-supermodular.  

Let~$\Eone := \Eone(\Pi)$ and~$\Etwo := \Etwo(\Pi, S)$.  For~$i=1,2$ we denote for a graph~$G$ by
\[
	f_i(G) := \begin{cases}
		1,& \text{ if } G\in \mathcal{E}_i \\
		0,& \text{ otherwise}
		\end{cases}
\]
the indicator function for the event~$\mathcal{E}_i$.  In the sequel we shall show that~$f_2$ is decreasing 
with respect to~$\le_\Pi$, and that~$f_1$ is increasing.  This will conclude the proof, as the above 
discussion yields that the conditions for (\ref{eq:myfkg}) are fulfilled.

First we prove that~$f_2$ is decreasing.  For this, it obviously suffices to show that if~$G\le_\Pi H$, 
then~$H\in \Etwo$ implies that~$G \in \Etwo$.  To see this, observe that~$H\in \Etwo$ implies that there is 
a set~$X\subset E(H;\,\Pi)$ such that
\[
		(E(H;\, \Pi)\setminus X)\cup S
		\text{ is triangle-free, and }
		|X| \le |S|.
\]
Due to~$G\le_\Pi H$, we have~$E(G;\,\Pi) \subseteq E(H;\,\Pi)$.  Let~$X' := X \cap E(G;\,\Pi)$, and note 
that~$|X'| \le |S|$.  Also, as~$(E(H;\, \Pi)\setminus X)\cup S$ is triangle free, so 
is~$(E(G;\, \Pi)\setminus X')\cup S$.  But this means~$G\in \Etwo$, as desired.

Finally, we prove that~$f_1$ is increasing. For this, we show that if~$G\le_\Pi H$, 
then~$G \in \Eone$ implies~$H \in \Eone$. Observe that by transitivity it is sufficient to consider the case 
that~$G$ and~$H$ differ in exactly one edge~$e$. The event~$G \in \Eone$ implies
\begin{equation}\label{eq:PiinEone}
	\gap(G;\;\Pi) \le r_0
	~\text{ and }~
	\forall\Pi' \text{ such that } \gap(G;\; \Pi,\Pi') \le r_0 : \dist(\Pi,\Pi') \le s_0.
\end{equation}
Now we make a case distinction. First assume that~$e$ joins two vertices in~$A$ or two vertices in~$B$.  
Then, due to (\ref{eq:ordering}), we have~$H = G\setminus\{e\}$.  Note that this implies that the size of 
a maximum bipartition satisfies~$b(G)-1\le b(H)\le b(G)$.  As furthermore~$E(G;\Pi)=E(H;\Pi)$, we thus 
have~$\gap(H;\, \Pi) \le \gap(G;\, \Pi) \le r_0$.  Now let~$\Pi'$ be a bipartition which has the 
property~$\gap(H;\, \Pi,\Pi') \le r_0$.  Observe that~$e(H;\, \Pi') \le e(G;\, \Pi')$.  We easily deduce 
\[
	\gap(G;\, \Pi,\Pi')
	\;=\; e(G;\,\Pi) - e(G;\,\Pi')
	\;\le\; e(H;\,\Pi) - e(H;\, \Pi')
	\;=\; \gap(H;\, \Pi,\Pi')
	\;\le\; r_0,
\]
which implies with (\ref{eq:PiinEone}) that~$\dist(\Pi,\Pi') \le s_0$.
	
Now assume that~$e$ joins a vertex in~$A$ with a vertex in~$B$. Then~$H = G\cup \{e\}$.  In this case we 
have~$b(G)\le b(H)\le b(G)+1$ and~$E(H;\,\Pi)=E(G;\,\Pi)\cup \{e\}$. This immediately implies
\[
	\gap(H;\, \Pi)
	= b(H) - e(H;\, \Pi)
	\le b(G) + 1 - (e(G;\, \Pi) + 1)
	= \gap(G;\, \Pi)
	\le r_0.
\]
Now let again~$\Pi'$ be a bipartition with~$\gap(H ;\; \Pi,\Pi') \le r_0$.  Note 
that~$e(G;\, \Pi') \ge e(H;\,\Pi') - 1$, as the edge~$e$ does not neccesarily join two vertices in 
different parts of~$\Pi'$.
%
This implies 
\[
	\op{gap}(G;\, \Pi,\Pi')
	= e(G;\,\Pi) - e(G;\, \Pi')
	\le (e(H;\, \Pi)-1) - (e(H;\, \Pi')-1)
	\le \gap(H ;\, \Pi,\Pi')
	\le r_0.
\]
Hence, (\ref{eq:PiinEone}) implies~$\dist(\Pi,\Pi') \le s_0$, as desired. This completes the proof.
\end{proof}

As a final ingredient in our proof we need estimates for the probabilities~$\Prob{\Gnp\in\Eone(\Pi)}$ 
and~$\Prob{\Gnp \in\Etwo(\Pi,S)}$. These are given by the next proposition.

\begin{proposition}\label{prop:finishingOff}
For all~$\Pi = (A,B)\in\BalPi_n$ and~$S\subseteq\binom{A}{2}\cup\binom{B}{2}$ with~$1\le |S|\le r_0$, and all 
sufficiently large~$n$, we have
\begin{equation}\label{eq:fOff}
	\Prob{\Gnp \in\Etwo(\Pi,S)} \le e^{-\frac{p^2n}{12}}
\end{equation}
and
\begin{equation}\label{eq:fOffII}
	\sum_{\Pi\in\BalPi_n}\Prob{\Gnp\in\Eone(\Pi)} \le n\binom{n}{s_0}.
\end{equation}
\end{proposition}

\begin{proof}
We first bound the probability of the event~$\Etwo$.  Fix an edge~$e\in S$ and note that, as~$\Pi$ is balanced, 
there exist at least
\[
	\min\{|A|, |B|\} \ge \frac{n}{2} - \frac{n}{100} = \frac{49}{100}n
\]
pairwise vertex-disjoint possible triangles across~$\Pi$ which contain the edge~$e$. 

Denote by~$Y$ the random variable which counts the number of those triangles in~$G_{n,p}$.  With the 
definition of~$\Etwo$ in (\ref{eq:E1E2}), we deduce
\[
	\Prob{G_{n,p} \in \Etwo(\Pi,S)} \le \Prob{Y \le |S|}.
\]
The probability that a triangle (with~$e$) is contained in~$G_{n,p}$ is~$p^2$; 
hence we have~$\E{Y} \ge \frac{49}{100}p^2n$.
On the other hand, observe that~$|S|\ll\E{Y}$ with our assumptions on~$p$.  A simple application of 
Lemma~\ref{len:cherny} yields \eqref{eq:fOff}, provided $n$ is sufficiently large.

Next we show (\ref{eq:fOffII}). Trivially, we have
\[
	\sum_{\Pi\in\BalPi_n}\Prob{\Gnp\in\Eone(\Pi)}
	=
	\sum_{\Pi\in\BalPi_n} \sum_{G \in \Eone(\Pi)}\Prob{\Gnp = G}
\]
Now we want to interchange the order of summation, such that the first sum goes over (a carefully chosen 
subset of) all graphs in~$\mathcal{G}_n$.  To achieve this, observe that the number of times the probability 
of a graph~$G$ is counted above is equal to the number of balanced partitions~$\Pi$ with the properties
\begin{equation}\label{eq:EoneNeg}
\begin{split}
	\gap(G;\;\Pi) \le &\; r_0, 
	\text{ and }
	\forall
	\Pi' \text{ such that } \gap(G;\; \Pi,\Pi') \le r_0 : \dist(\Pi,\Pi') \le s_0.
\end{split}
\end{equation}
In the following we argue that we can construct all such partitions~$\Pi$ by taking the canonical optimal 
partition~$\Pi^*$, and modifying the parts of~$\Pi^*$ in at most~$s_0$ vertices.  To see this, let $\Pi$ 
be any partition~$\Pi$ fulfilling (\ref{eq:EoneNeg}); then we have~$\gap(\Pi, \Pi^*) \le r_0$, 
which implies~$\dist(\Pi, \Pi^*)\le s_0$. 

Hence, as there are precisely~$\sum_{t\le s_0}\binom{n}{t}$ ways to choose at most~$s_0$ vertices, which 
change the class they belong to, and for sufficiently large~$n$ the inequality~$s_0 \le \frac{n}{2}$ holds, 
we can conclude
\[
	\sum_{\Pi\in\BalPi_n}\Prob{\Gnp\in\Eone(\Pi)}
	\le
	\sum_{G\in\mathcal{G}_n}
		\Big(
			\Prob{\Gnp = G} \cdot 
			\sum_{\substack{\Pi\in\BalPi_n \\ \Pi \text{ fulfills (\ref{eq:EoneNeg})}}} 1
		\Big)
	\le
	n\binom{n}{s_0}.
\]
\end{proof}

\begin{proof}[Proof of Theorem \ref{thm:extremalMain}]
Recall that~$n^{-1/250} \le p \le \frac12$, which implies the bounds~$r_0 \le n^{13/250}$ 
and~$s_0 \le n^{9/10}$ for sufficiently large $n$.  As there are at most~$\binom{n^2}{|S|}$ ways to choose a 
set~$S$ of edges out of all possible edges, the proof of the theorem can be completed with 
inequality~\eqref{eq:getRidOfBadEvents} and Propositions \ref{prop:tranferToE1E2}, 
\ref{prop:FKGApplication}, and \ref{prop:finishingOff} as follows:
\[
\begin{split}
	\Prob{\exists T\in\mathcal{T}(\Gnp): T \text{ not bipartite}}& \\
	&\le
	\Prob{\mathcal{B}} + \sum_{\substack{\Pi\in\BalPi_{n}\\ S:|S|\le r_0}}
			\Prob{\Gnp\in\Eone(\Pi)} \cdot \Prob{\Gnp\in\Etwo(\Pi,S)} \\
	&\hspace{-1cm}\le
	o(1) + e^{-\frac{p^2n}{12}} \cdot \binom{n^2}{r_0} \cdot 
	\sum_{\Pi\in\BalPi_{n}}\Prob{\Gnp\in\Eone(\Pi)} \\
	&\hspace{-1cm}\le
	o(1) + \exp
		\Big\{
		- \frac{p^2n}{12}
		+ 2r_0\log n
		+ (s_0 + 1)\log n
		\Big\}
	= o(1),
\end{split}
\]
whenever~$n$ is chosen sufficiently large.
\end{proof}


\section{Larger complete graphs}\label{sec:generalizations}

Let~$\ell\ge 2$. Here we show how the proofs of the previous sections can be adapted in order to prove Theorem \ref{thm:cliques}. As a tool, we will use a ``higher-dimensional'' variant of Theorem~\ref{thm:propertiesOfCuts}, see Theorem \ref{thm:propertiesOfellCuts} below. Before we state it, we need to define the notion of distance for two~$\ell$-partitions, which is a straightforward generalization of the case~$\ell = 2$:
\[
	\dist(\Pi, \Pi')
	:= \min_{\substack{\pi:[\ell]\rightarrow[\ell] \\ \pi\text{ bijection}}}
			\sum_{\substack{1\le i\le\ell \\j: \pi(j) \not= i}} \card{V_i \cap V'_{\pi(j)}}.
\]
The notion of the gap of two~$\ell$-partitions is defined in the obvious way. The following theorem is a statement about the structure of the set of (near-)optimal~$\ell$-partitions of the uniform random graph.
\begin{theorem}\label{thm:propertiesOfellCuts}
Let~$\ell \ge 2$. There are constants~$C  = C(\ell) > 0$ and~$\eps_0 = \eps_0(\ell)> 0$ such that the following holds for sufficiently large~$n$. Let~$0 \le \eps \le \eps_0$,~$n^{-1/\ell} \ll p \le \frac12$ and~$M = p\binom{n}{2}$. Furthermore, let~$r\ge 1$ satisfy~$r \ll (pn)^{1/8}$ and~$\omega \gg 1$, and define
\[
	s_0 := C \cdot \omega \cdot r^4 \cdot \sqrt{np^{-1}}.
\]
Then
\[
	\Prob{\exists \Pi : \gap(G_{n,M}; \; \Pi) = r-1 \;\text{ and }\; \dist(G_{n,M}; \; \Pi) \ge s_0}
	~\le~
	\omega^{-1}.
\]
\end{theorem}
The proof of Theorem \ref{thm:propertiesOfCuts} can easily be adapted to show the above theorem. The sole difference is that instead of considering bipartitions we have to consider~$\ell$-partitions. In this case it is readily seen that~\eqref{eq:expectedChangeLowerBound} is still valid, i.e., the expected increase of the size of a maximum~$\ell$-partition can be estimated from below by
\[
	\E{b(\GnMt_i) - b(\GnM)} \ge \E{\mathcal{I}_{t_i}} + \E{\io},
\]
where~$b$ denotes the size of maximum~$\ell$-partition,~$\mathcal{I}_{t_i}$ the increase of the size of a fixed maximum~$\ell$-partition~$\Pi^*$ of~$\GnM$, and~$\io$ the indicator variable for the event that a partition~$\Pi$ with gap~$r-1$ becomes an optimal bipartition after adding~$t_i$ random edges to~$\GnM$.~$\E{\mathcal{I}_{t_i}}$ can be routinely estimated from below. Moreover, to estimate~$\E{\io}$ we proceed exactly as in the proof of Theorem~\ref{thm:propertiesOfCuts}: we estimate the probability that the number of edges added that increased the size of~$\Pi$ is at least~$r$ plus the number of edges which increased the size of~$\Pi^*$. We leave the solely technical but straightforward details to the reader.

With the above observations, it is not very surprising that the proof of Theorem \ref{thm:extremalMain} can be adapted in order to prove Theorem \ref{thm:cliques} -- in fact, it turns out that it does not make a difference for our proofs if we consider~$\ell$-partitions instead of bipartitions of the~$G_{n,p}$. However, some details are significantly more tedious than it is above the case, and we shall elaborate more on this issue.

We proceed in three steps, mimicking the proof Theorem \ref{thm:extremalMain}. Let~$\mathcal{F}(\Gnp)$ denote the set of the maximum~$\mathcal{K}_\ell$-free subgraphs of~$\Gnp$, and let~$F\in\mathcal{F}(\Gnp)$. First, we prove a statement similar to that of Lemma \ref{lem:pn2}: for every~$\eps > 0$ we can find a.a.s.\ a partition~$\Pi = (V_1,\dots,V_{\ell-1})$ such that all but~$\eps pn^2$ edges of~$F$ go across~$\Pi$, and all parts of~$\Pi$ have approximately the same size. The proof is, as in the case of maximum triangle-free graphs, an application of the sparse version of Szemeredi's regularity lemma (Theorem \ref{thm:szem_sparse}) and the probabilistic embedding lemma (Theorem \ref{thm:relaxedKLR}), followed by an application of the stability lemma which results in the desired~$(\ell-1)$-partition. The calculations differ only in minor technical details, which are again left to the reader.

Second, we show that in fact we can find a better~$(\ell-1)$-partition, i.e., we prove a general version of Lemma~\ref{lem:almostlpartite}.
\begin{lemma}\label{lem:almostell-1partite}
For~$\ell \ge 4$ and~$p\ge n^{-\frac{1}{100\ell^3}}$ the following holds a.a.s. For every~$F\in\mathcal{F}(G_{n,p})$ there is a partition~$\Pi_F= \Pi = (V_1, \dots, V_{\ell-1})$ of the vertex set such that all but~$2p^{-5\ell^2}\log^2 n$ edges of~$F$ go across~$\Pi$. Furthermore,~$|V_i| = \frac{n}{\ell-1} + o(n)$ for~$1 \le i\le \ell-1$.
\end{lemma}
Before we give the proof details let us state an auxiliary result that we will use several times. It is a statement about the number of pairwise edge-disjoint copies of
complete graphs in subgraphs of the~$\Gnp$.
\begin{proposition}\label{prop:KellhelperII}
$G_{n,p}$ has for every~$t \ge 3$ a.a.s. the following property. There are constants~$c = c(t)$,~$c' = c'(t)$ such that for every~$t$ disjoint subsets of its vertices~$S_1, \dots, S_t$, with~$|S_1| := s_1 \ge cp^{-t^2}\log n$ and for~$i\ge 2$~$|S_i| =: s \ge |V_1|$, the number of pairwise edge-disjoint~$K_t$'s with one vertex in each~$S_i$ is at least~$c'ps_1s$ and at most~$\le 2ps_1s$.
\end{proposition}
\begin{proof}
The upper bound is easy to obtain, as the number of edges between the sets~$S_1$ and~$S_2$ is a.a.s.\ at most~$2ps_1s$, and each copy of $\mathcal{K}_t$ contains one of those edges. For the lower bound, let us fix~$S_1, \dots, S_t$. We apply Lemma~\ref{lem:randomnessGraham} recursively several times with~$c = 1/2$. First for $U\leftarrow S_t$ and $k\leftarrow t-1$ (to obtain an exceptional set~$X_{t-1}$), then for the set $U\leftarrow S_{t-1}\setminus X_{t-1}$ and $k\leftarrow t-2$ (to obtain an exceptional set~$X_{t-2}$), and so on, until $U\leftarrow S_2\setminus(X_{t-1}\cup\ldots\cup X_2)$ and $k\leftarrow 1$ (to obtain an exceptional set~$X_{1}$). Let ~$S_i' := S_i\setminus(\cup_{j=i}^{t-1}X_{j})$, and note that for sufficiently large~$n$ we have~$|S_i'| \ge \frac12|S_i|$, as we have $|X_i| =\mathcal{O}(p^{-i}\log n)$. 
Note also that by construction any tuple $(v_1,\dots,v_j)$ with $v_i\in S_i'$, for $1\le i\le j$, has at least~$\frac12p^j|S_j'|$ common neighbors in~$S_{j+1}$.
Hence, we deduce that there exist at least
\[
 |S_1'| \cdot  \frac12p |S_2'| \cdot   \frac12p^2 |S_3'| \cdot \ldots \cdot  \frac12p^{t-1} |S_t'| \ge 2^{-2t} p^{\binom{t}2} \prod_{i=1}^t |S_i| \ge 2^{-2t} p^{\binom{t}2} s_1s_2^{t-1}
\]
$\mathcal{K}_t$'s with one vertex in each~$S_i'$. To transfer this lower bound on the {\it total number} of $\mathcal{K}_t$'s into a lower bound on the maximum size of
pairwise edge-disjoint $\mathcal{K}_t$'s, we now establish an upper bound on the number of $\mathcal{K}_t$'s that contain a fixed edge $e$.            

Let~$e=\{v_{i_1},v_{i_2}\}$ be an edge that joins a vertex in~$S_{i_1}'$ to one in~$S_{i_2}'$. The number of triangles with~$e$ and a vertex~$v_{i_3}\in S_{i_3}$, $i_3\not\in\{i_1,i_2\}$, is at most~$\frac32p^2|S_{i_3}'|$, as the common neighborhood of the endpoints of~$e$ in~$S'_{i_3}$ has at most this size. Repeating this argument for the vertices $v_{i_1}, \ldots, v_{i_{k-1}}$ and an index $i_{k}\not\in \{i_1,\ldots,i_{k-1}\}$ we see that there exist at most~$\frac32p^{k-1}|S_{i_{k}}'|$ vertices $v_{k}$ in $S_{i_{k}}'$ that form a $\mathcal{K}_k$ with  $v_{i_1}, \ldots, , v_{i_{k-1}}$. Hence, the number of  $\mathcal{K}_t$'s that contain the edge~$e$
is bounded from above by
\[
\left(\frac32\right)^{t-2} \cdot\prod_{i=2}^{t-1}p^i\cdot \prod_{{1\le i\le t}\atop{i\not\in\{i_1,i_2\}}}^t |S_i| \le \left(\frac32\right)^{t-2}p^{\binom{t}2-1} s_2^{t-2},
\]
which concludes the proof of the lemma in a straightforward way. 
\end{proof} 

\begin{proof}[Proof of Lemma~\ref{lem:almostell-1partite}]
The main strategy is very similar to the strategy used in the proof of Lemma~\ref{lem:almostlpartite}, but we have to make some important modifications. We will sketch only the relevant steps, and the missing details can easily be filled in by considering the respective steps in the original proof. 

First of all, note that due to the discussion before Lemma~\ref{lem:almostell-1partite} there is a.a.s.\ a partition~$(V_1^F, \dots, V_{\ell-1}^F)$ of~$F$ such that
\[
	\sum_{i=1}^{\ell-1} e(F;\, V_i^F) \le {\eps^5} pn^2, \text{ and } |V_i^F| = \frac{n}{\ell-1}\pm \eps^5n,
\]
for any sufficiently small positive~$\eps$. Moreover, with similar arguments as in the proof of Lemma~\ref{lem:almostlpartite} it can be shown that every~$(\ell-1)$-partition~$(P_1, \dots, P_{\ell-1})$ of~$\Gnp$ with gap at most~$\eps^5pn^2$ satisfies a.a.s.~$(1 - \eps)\frac{n}{\ell-1}\le |P_i| \le (1 + \eps)\frac{n}{\ell-1}$. We will assume that~$\Gnp$ has all those properties simultaneously, and make also all other additional assumptions on~$G_{n,p}$ made in Lemma~\ref{lem:almostlpartite}.

In accordance with Lemma~\ref{lem:almostlpartite}, we call a partition~$\Pi$ \emph{optimal} with respect to~$F$, if~$e(F;~\Pi)$ attains its maximum over all possible partitions. By exploiting our assumptions we see that all optimal~$(\ell-1)$-partitions~$\Pi = (V_1, \dots, V_{\ell-1})$ of~$F$ satisfy~$(1 - \eps)\frac{n}{\ell-1} \le |V_i|\le (1 + \eps)\frac{n}{\ell-1}$ and~$\sum_{i=1}^{\ell-1} e(F;\, V_i) \le {\eps^5} pn^2$.

Let~$\Pi = (V_1, \dots, V_{\ell-1})$ be any fixed optimal~$(\ell-1)$-partition of~$F$. A \emph{horizontal} edge is an edge in~$F$ that joins two vertices in the same~$V_i$, and a \emph{missing} edge joins in~$G_{n,p}$ two vertices in different~$V_i$'s, but is not contained in~$F$. As in Lemma~\ref{lem:almostlpartite} we can then define the horizontal degree and missing degree of any vertex~$v$ with respect to~$F$.

Next, for all~$2 \le i\le \ell-1$ we denote by~$\mathcal{B}_i(U)$ a (minimum) set of vertices, such that for every~$v_1, \dots, v_i \not\in\mathcal{B}_i(U)$ we have~$\big||\bigcap_{j=1}^{i} \Gamma(G_{n,p};\, v_j)\cap U| - p^i|U|\big| \le \frac14 p^i|U|$. Let
\[
	X_1^{V_i}
	=
	\left(\bigcup_{j=1}^{\ell-1} \mathcal{B}_1(V_j) \cup \dots \cup \mathcal{B}_{\ell-1}(V_j)\right)\;\cap\; V_i,
	\quad ~ ~ \text{ and } X_2^{V_i}, X_3^{V_i} \text{ remain as in Lemma~\ref{lem:almostlpartite}}.
\]
Moreover, define~$V_i^0 := V_i$, and~$V_i^{j} := V_i^{j-1} \setminus X_j^{V_i}$, and~$X_i := \cup_{k=1}^{\ell-1}X_i^{V_k}$. Note that due to Lemma~\ref{lem:randomnessGraham} we have~$|X_1| = \mathcal{O}(p^{-\ell+1})$, and that the number of missing edges
in~$F$ incident to at least one vertex in~$|X_3|$ is at least~$\frac52|X_3|\eps pn$. From this we readily obtain that~$|X_3| \le \eps n$, as the number of missing edges is at most the number of horizontal edges: otherwise~$F$ would not be a maximum~$\mathcal{K}_\ell$-free graph.

We now proceed with the following steps, mimicking the proof of Lemma~\ref{lem:almostlpartite}.
\begin{itemize*}
	\item[(i)] We first show that~$X_2$ is small, i.e.,~$|X_2| \le \eps p^{-2}$.
	\item[(ii)] Set~$H_3 := \cup_{i=1}^{\ell-1} E(F;\, V_i^3)$. We show that~$\card{H_3} \le p^{-\ell^2}n\log n$.
	\item[(iii)] We show~$|X_3| \le p^{-\ell^2-2}\log n$.
	\item[(iv)] We use (ii) to show that for all~$v \in \cup_{i=1}^{\ell-1}V_i^3$ we have~$d_H(v) \le p^{-2\ell^2}\log n$.
	\item[(v)] Then we show that in fact~$|H_3| \le p^{-5\ell^2}\log^2n$.
	\item[(vi)] Finally, we show~$\sum_{i=1}^{\ell-1}e(F;\, V_i) \le 2p^{-5\ell^2}\log^2n$.
\end{itemize*}
The proofs of these statements follow along similar lines as the proofs for the corresponding statements in Lemma~\ref{lem:almostlpartite}, they just get technically more envolved. To see (i), suppose that~$|X_2^{V_1}| \ge t_0 = c_\ell(\frac\eps{p})^2$, where~$c_\ell$ will be specified later, and fix a set $X`\subseteq X_2^{V_1}$ of size $|X'|=t_0$. The number of missing edges at a vertex $v\in X'$ can be bounded from below by the maximum number of pairwise edge-disjoint~$K_\ell$'s, which have one endpoint in each of~$V_2\cap \Gamma(G_{n,p};\, v), \dots, V_{\ell-1}\cap \Gamma(G_{n,p};\, v)$, and two endpoints in~$V_1$, namely~$v$, and one of its neighbors in~$V_1$. Observe that~$d_H(v) \le d(F; v, V_i)$ for all~$v\in X'$ and all~$2\le i\le \ell-1$, i.e., the neighborhoods of~$v$ in all~$V_i$ have size at least~$\eps p n$. Let
\[
\begin{split}
\mathcal{F}(v) &:= \text{ maximum set of pairwise edge-disjoint~$\mathcal{K}_\ell$'s in~$G_{n,p}$, which} \\
&\text{contain~$v$ and one of its neighbors in each~$(V_j)_{j=1\dots \ell-1}$},
\end{split}
\]
and set~$m_0 = |\bigcup_{v\in X_2^{V_1}} \mathcal{F}(v)|$. Note that the set $\mathcal{F}(v)$ may not be uniquely defined; this does not matter, we just fix maximum sets arbitrarily. 

Observe that~$\mathcal{F}(v)$ is equal to a maximum set of pairwise edge-disjoint $\mathcal{K}_{\ell-1}$'s between the neighborhoods of~$v$. As these neighborhoods are all of size at least~$\eps pn$, by applying Proposition~\ref{prop:KellhelperII} we see that a.a.s.\ we have for all $v\in X'$ that
$|\mathcal{F}(v)| \ge C_\ell p\cdot(\eps pn)^2$, for some~$C_\ell > 0$. Now we perform a similar calculation as in \eqref{eq:mr}.
Let~$m_0 := |\bigcup_{v\in X'}\mathcal{F}(v_i)|$ and observe that the number of missing edges is at least~$m_0$. We apply the inclusion-exclusion principle:
\[
	m_0\;
	\ge
	\sum_{v\in X'}  \card{\mathcal{F}(v)} - \sum_{v,w\in X', v\not=w} \card{\mathcal{F}(v) \cap \mathcal{F}(w)}
\]
where~$\mathcal{F}(v) \cap \mathcal{F}(w)$ denotes the set of {\sl edges} that are contained in the intersection of copies  of $\mathcal{K}_\ell$'s in $\mathcal{F}(v)$ and $\mathcal{F}(w)$. Observe that~$\card{\mathcal{F}(v) \cap \mathcal{F}(w)}$ can be crudely bounded from above by the number of edges between the common neighborhoods of~$v$ and~$w$ in~$V_2, \dots, V_{\ell+1}$. But these neighborhoods have size at most~$2p^2n$ (as a.a.s.\ no two vertices have a larger common neighborhood in~$\Gnp$), and the number of edges between any two sets of vertices of at most this size in~$G_{n,p}$ is a.a.s.~$\le 2p(2p^2n)^2$, which implies~$\card{\mathcal{F}(v) \cap \mathcal{F}(w)} \le 4\ell^2p(p^2n)^2$. By choosing~$c_\ell = \frac{C_\ell}{8\ell^2}$ we see that the assumption $|X'|=t_0 = c_l \eps^2 p^{-2}$ implies that~$m_0 \ge \frac1{2} c_\ell C_\ell\eps^4pn^2$, which contradicts the fact~$m_0 \le \eps^5pn^2$, whenever~$\eps$ is sufficiently small.

Next we show how (ii) can be proved. Assume that the number of horizontal edges is maximized in~$V^3_1$. We adapt the definition of a ``chord'': here, a chord consists of~$\ell$ vertices~$\mathcal{V} = \{x, y, v_2, \dots , v_{\ell-1}\}$, such that~$\{x, y\}$ is a horizontal edge in~$V^3_1$, and~$v_i \in V_i^3$; additionally, one of the edges joining vertices in~$\cal V$ (except for~$\{x,y\}$) is a missing edge, i.e., it is contained in~$\Gnp$ but not in~$F$, and all other edges are in~$\Gnp$. See Figure~\ref{fig:l_chord} for an illustration.

\begin{figure}
\centering
\includegraphics[height=31mm]{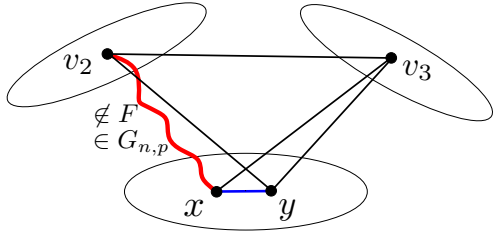}
\caption{A chord in the case~$\ell=4$. The black edges are in~$G_{n,p}$, and the red edge is missing.}
\label{fig:l_chord}
\end{figure}

Our objective is to derive an upper and a lower bound for the number~$N$ of chords, which will immediately imply the bound on~$|H_3|$ claimed in (ii). We begin with the lower bound. Note that we may assume that
$|H_3| > p^{-\ell^2}n\log n$, as otherwise there is nothing to prove. By symmetry we may in addition assume that $|E(F,V_1^3)|\ge |H_3| / \ell$.
 Due to Lemma~\ref{lem:randomness} there is a set $Z$ of $\mathcal{O}(p^{-2})$ vertices such that every pair in $V_1^3\setminus Z$ has at least at least~$\frac{3p^2|V_2^3|}{4} \ge \frac{3p^2n}{4\ell}$ common neighbors in $V_2^3$. As at most $|Z|\cdot |V_1^3|=\mathcal{O}(p^{-2}n)$ pairs of vertices in $V_1^3$ contain a vertex from $Z$, we deduce that the number of triangles in~$\Gnp$ connecting any edge $\{x,y\}$ in $V_1^3$ with a vertex $v_2\in V_2^3$ is at least
\begin{equation}\label{eq:trianglesLowerBound}
	\left(
		\Big(\frac{|H_3|}\ell - \mathcal{O}(p^{-2}n)\Big)
	\right)
	\frac{3p^2n}{4\ell}
	\stackrel{(|H_3| > p^{-\ell^2}n\log n)}{\ge}
	|H_3| \cdot \frac{p^2n}{2\ell^2}.
\end{equation}
Similarly, the number of~$\mathcal{K}_4$'s connecting~$x,y,v_2$ with a~$v_3\in V_3^3$ is at least~$\frac{3p^3n}{4\ell}$, except for at most~$\mathcal{O}(p^{-3}n^2)$ triples~$x,y,v_2$. Exploiting \eqref{eq:trianglesLowerBound} we obtain that the number of~$\mathcal{K}_4$'s with an edge in $V_1^3$, a vertex in $V_2^3$, and a vertex in $V_3^3$ is at least
\[
	 \left(
		\Big(|H_3| \cdot \frac{p^2n}{2\ell^2} - \mathcal{O}(p^{-3}n^2)\Big)
	\right)
	\frac{3p^3n}{4\ell}
	\ge
	|H_3| \cdot \frac{p^5n}{2^2\ell^3}.
\]
This process can be continued to count~$\mathcal{K}_5$'s,~$\mathcal{K}_6$'s and so on. 
We obtain with room to spare that the number of~$\mathcal{K}_\ell$'s is at least
\begin{equation}\label{eq:ellChordsLowerBound}
|H_3| \cdot \frac{n^{\ell-2} p^{\binom{\ell}{2}-1}}{(2\ell)^{\ell}},
\end{equation}
which provides us with the desired lower bound for the number of chords. 

To obtain an upper bound for~$N$ we partition the set of chords into three classes and derive upper bounds for each one separately.
\begin{itemize}
	\item[(A)] Chords where~$d_{H_3}(x) \le p^{-2}\log n$, and there is a~$2 \le i \le \ell-1$ such that~$\{x,v_i\}$ is missing. We count those as follows. Let~$x\in V_1^3$ and~$y \in \Gamma(F;\,x) \cap V_1^3$. Then an upper bound for the number of chords with~$x$ and~$y$ is given by~$\ell^2$ multiplied with the number of~$\mathcal{K}_\ell$'s in~$G_{n,p}$, that have one vertex in each~$V_i$ ($2\le i\le V_{\ell-1}$) and $x,y$. To estimate this number we use a similar argument as above: the number of~$\mathcal{K}_3$'s with~$x,y$ and a vertex~$v_2\in V_2$ is at most~$2p^2n$, as no two vertices have a larger neighborhood in~$G_{n,p}$. Similarly, the number of~$\mathcal{K}_4$'s with~$x,y,v_2$ and~$v_3\in V_3$ is at most~$2p^3n$, and so on. Putting everything together yields that the number of (A)-chords is at most
	\[
	\begin{split}
		&\sum_{\substack{x\in V_1^3 \\ d_{H_3}(x) \le p^{-2}\log n}}
			\sum_{y \in \Gamma(F;\,x) \cap V_1^3}
				(2p^2n)\cdot(2p^3n)\cdots(2p^{\ell-1}n) \\
		&\le |V_1^3|\cdot p^{-2}\log n\cdot (2n)^{\ell-2}p^{\binom{\ell}{2}-1}
		~\le~ (2n)^{\ell-1} p^{\binom{\ell}{2}-3} \log n .
	\end{split}
	\]
	\item[(B)] Chords where~$d_{H_3}(x) \ge p^{-2}\log n$, and there is a~$2 \le i \le \ell-1$ such that~$\{x,v_i\}$ is missing. Let~$S = \Gamma(G_{n,p};\, x) \setminus (V_1^3 \cup \Gamma(F;\, x))$ and note that~$d_M(x) = |S|$. Due to our assumptions~$G_{n,p}$ is such that except of at most~$Cp^{-1}\log n$ vertices, every vertex in~$S$ has at most~$\frac54pd_{H_3}(x)$ neighbors in~$\Gamma(F;\, x)\cap V_1^3$. Hence, the number of~$\mathcal{K}_3$'s with~$x$, a neighbor of~$x$ in~$V_1^3$ and a vertex adjacent to a missing edge at~$x$ is at most
	\[
		d_M(x) \cdot \frac54p\, d_{H_3}(x) + Cp^{-1}\log n \cdot d_{H_3}(x) 
		\stackrel{(d_M(x) \le 6\eps pn)}{\le} 8\eps p^2 n \cdot d_{H_3}(x).
	\]
	With a similar argument as in (A) we easily see that the number of~$\mathcal{K}_4$'s with~$x$, a neighbor of~$x$ in~$V_1^3$, a vertex adjacent to a missing edge at~$x$, and a vertex in any of~$V_i^3$'s (different from the ones where a vertex is already taken from) is at most~$2p^3n$, as this is the maximum number of common neighbors of any three vertices in~$G_{n,p}$. In the same way we can count~$\mathcal{K}_5$'s,~$\mathcal{K}_6$'s, etc. To conclude, the number of (B)-chords is at most
	\[
		\sum_{\substack{x\in V_1^3 \\ d_{H_3}(x) > p^{-2}\log n}}
		8\eps p^2 n \cdot d_{H_3}(x) \cdot (2p^3n)\cdots (2p^{\ell-1}n)
		~\le~
		8\eps(2n)^{\ell-2} p^{\binom{\ell}{2}-1} \cdot |H_3|.
	\]
	\item[(C)] Chords with the property that there are indexes~$2 \le i,j \le \ell-1$ such that~$\{v_i, v_j\}$ is missing. We count (C)-chords as follows. 
	In a first step we count $\mathcal{K}_4$'s with vertices $x,y,v,v'$ such that $x,y$ form a horizontal edge in $V_1^3$, $v,v'$ come from different sets $V_i^3$, $i\ge 2$ such that $\{v, v'\}$ form a missing edge. In a second step we then bound the number of ways in which these $\mathcal{K}_4$ can be extended to a 
	 $\mathcal{K}_\ell$.
	 
	Let~$e = \{x,y\}\in H_3$. Note that there are at most~$2p^2n$ candidates for~$v$, as $v$ needs to be connected to $x$ and $y$. Let~$M_v$ be the set of missing neighbors of $v$, i.e., the set of vertices incident to missing edges $e = \{v,w\}$, where $w\in V_i^3$ for some $i\ge 2$.  We now distinguish various cases. Case (a): $M_v$ is small, more precisely, $|M_v|< C p^{-2}\log n$, where $C$ is an appropriately chosen constant, see below. Then we can obviously bound the number of choices for $v'$ by $|M_v| = \mathcal{O}(p^{-2}\log n)$. 
	Case (b): $|M_v| \ge C p^{-2}\log n$ and $x,y$ both do not belong to the exceptional set $Q_{M_v}$ defined in Lemma~\ref{lem:randomnessGraham} for $M_v$, where $C=C(2,1/2)$. Then Lemma~\ref{lem:randomnessGraham} implies that there are at most $\frac32p^2 |M_v|\le6\eps p^3n$ ways to choose~$v'$, as $d_M(v)\le 6\eps pn$.
	Finally, we treat in case (c) those tuples $x,y,v$ where at least one of $x$ and $y$ are contained in the exceptional set with respect to $M_v$. Note that
	as Lemma~\ref{lem:randomnessGraham} implies that the size of the exceptional set is bounded by $\mathcal{O}(p^{-2}\log n)$, the number of such
	triples is thus at most   $\mathcal{O}(n^2p^{-2}\log n)$. And for each such triple we have at most $2pn$ ways to choose a neighbor $v'$. Combining these numbers we see that there are at most
	\begin{eqnarray*}
	|H_3| \cdot p^2n\cdot  \mathcal{O}(p^{-2}\log n) + |H_3| \cdot p^2n\cdot  6\eps p^3n + \mathcal{O}(n^2p^{-2}\log n) \cdot 2pn&&\\
	\le\;|H_3| \cdot  8\eps p^5n^2 +  \mathcal{O}(p^{-1}n^3\log n)&&
	\end{eqnarray*}
	ways to chose~$x,y,v,v'$, if $n$ is sufficiently large.
	Having chosen~$x,y,v,v'$, the number of~$\mathcal{K}_\ell$'s containing those vertices, and one vertex in each of the remaining $V_i^3$'s is at most~$(2n)^{\ell-4}p^{\binom\ell2 - 6}$; this is seen by exactly the same counting argument that we have already used in (A) and (B). Putting everything together yields that the number of (C)-chords is at most
	\[
	|H_3| \cdot  2\eps p^{\binom\ell2 - 1}(2n)^{\ell-2} +  \mathcal{O}(p^{\binom\ell2 - 7}n^{\ell-1}\log n).
	\] 
\end{itemize}
	By combining the results from (A), (B) and (C) with \eqref{eq:ellChordsLowerBound} we see that~$|H_3|$ satisfies
\[
	|H_3| \cdot \frac{n^{\ell-2} p^{\binom{\ell}{2}-1}}{(2\ell)^{\ell}}
	\le 
	|H_3| \cdot 10\eps (2n)^{\ell-2} p^{\binom{\ell}{2}-1}
	+ \mathcal{O}(p^{\binom\ell2 - 7}n^{\ell-1}\log n),
\]
from which (ii) follows readily for large~$n$ (recall that $\ell \ge 4$).

The proof of (iii) is essentially identical to the proof of the analogous statement in Lemma~\ref{lem:almostlpartite}, we just have to adapt the powers of $p$. We omit a detailed exposition. To see (iv), observe that the total number of horizontal edges in~$F$ is due to (i)-(iii) at most
\[
	|H| \le |H_3| + |X_1 \cup X_2\cup X_3|n < 2p^{-\ell^2-2}n\log n.
\]
Suppose that there is a vertex~$v$ in~$V_1^3$ with~$d_{H_3}(v) \ge p^{-2\ell^2}\log n$; by symmetry vertices in $V_i^3$ for $i\ge 2$ can be handled analogously. We will show that this implies that the number of missing edges is at least~$2p^{-\ell^2-2}n\log n$, which contradicts the bound on~$|H|$ derived above.

In order to give a lower bound for the number of missing edges, we estimate from below the maximum number of edge-disjoint~$\mathcal{K}_\ell$'s in~$G_{n,p}$, which contain~$v$, one of the vertices counted in~$d_{H_3}(v)$, and~$\ell-1$ vertices in~$\Gamma(F;\, v, V\setminus V_1)$, such that there is precisely one vertex in each~$V_2, \dots, V_{\ell-1}$. For this we count the maximum number of edge-disjoint~$\mathcal{K}_{\ell-1}$'s between the sets of vertices~$\Gamma(F;\, v, V_2), \dots \Gamma(F;\, v, V_{\ell-1})$, and~$\Gamma(F;\, v, V_1^3)$ in~$G_{n,p}$. 

Note that the definition of the exceptional sets $X_j^{V_1}$ imply that we have~$d(F;\, v, V_i) \ge c_\ell pn$, for all $i\ge 2$, where~$c_\ell>0$ is an appropriate constant depending only on~$\ell$. Furthermore, recall that~$|\Gamma(F;\, v, V_1^3)|= d_{H_3}(v) \ge p^{-2\ell^2}\log n$ due to our assumption that $v$ is a vertex with such a high horizontal degree. We apply Proposition \ref{prop:KellhelperII} with~$k=\ell-1$,~$S_1 = \Gamma(F; v, V_1^3)$, and~$S_i = \Gamma(F; v, V_i)$ (truncated to their first~$c_\ell pn$ vertices) for $i\ge 2$, which yields that there is a constant~$c>0$ such that there are at least~$cp^2d_{H_3}(v)n$ pairwise edge-disjoint copies of~$\mathcal{K}_{\ell-1}$ with one endpoint in~$\Gamma(F;\, v, V_1^3)$ and in each~$\Gamma(F;\, v, V_i)$. But then the number of missing edges is~$\ge cp^{-2\ell^2+2}n\log n$, which completes the proof of (iv).

Next we prove (v).
Let~$m$ be the number of missing edges in~$F$. Our aim is to show that~$m \ge |H_3| \cdot \frac{p^{3\ell^2}}{\log n}n$, and hence~$|H_3|$ must satisfy
\[
	|H_3| \cdot \frac{p^{3\ell^2}}{\log n}n \le 2p^{-\ell^2-2}n\log n, \qquad \textrm{as } |H|< 2p^{-\ell^2-2}n\log n.
\]
This completes the proof of (v). To show the claimed bound for~$m$, assume that the number of edges in $H_3$ is maximized in $V_1^3$, and let~$R \subseteq H_3$ be a \emph{matching} of maximum cardinality that joins vertices from $V_1^3$. Note that by using (iv) we obtain~$|R| \ge \frac1\ell|H_3|p^{2\ell^2}(\log n)^{-1}$. We now proceed in two steps. First, we bound from below the number of~$\mathcal{K}_\ell$'s in~$G_{n,p}$, which contain an edge in~$R$, and a vertex in each of the sets~$V_2^3, \dots, V_{\ell-1}^3$. In the second step, we estimate from above the maximum number of~$\mathcal{K}_\ell$'s, that contain any edge in~$R$, and an additional (fixed) edge~$e'$ connecting two vertices in the sets~$V_i^3$ and~$V_j^3$, where~$1 \le i < j\le \ell-1$. By dividing these two numbers (and by dividing the result by~$\binom{\ell}2$) we readily obtain a lower bound for the number of missing edges. 

To obtain the first goal note that the number~$\tau_3$ of triangles with~$e\in R$ and a vertex~$v_2\in V_2^3$ is at least~$\frac34p^2|V_2| - |X_1\cup X_2 \cup X_3| \ge c_3p^2n$, for some~$c_3 > 0$. Similarly, the number~$\tau_i$ of~$\mathcal{K}_i$'s with~$e$ and~$i-2$ vertices~$v_2\in V_2^3$, \dots,~$v_{i-1}\in V_{i-1}^3$ is at least
\[
	\tau_i
	\ge \tau_{i-1} \cdot \left(\frac34 p^{i-1} \cdot |V_{i-1}| - |X_1\cup X_2 \cup X_3|\right)
	\ge \dots
	\ge c_i \cdot n^{i-2}p^{\binom{i}2-1},
\]
where~$c_i$ depends only on~$c_3$ and on~$i$. Setting~$i=\ell$ yields that the number of~$\mathcal{K}_\ell$'s in~$G_{n,p}$, which contain an edge in~$R$, and a vertex in each of the sets~$V_2^3, \dots, V_{\ell-1}^3$, is at least~$|R|c_\ell p^{\binom\ell2-1}n^{\ell-2} \ge |H_3| \cdot p^{3\ell^2-1}\cdot (\log n)^{-1} \cdot n^{\ell-2}$.

To obtain the second goal we distinguish two cases: either~$e'$ has a common endpoint with one edge in~$R$, or~$i\ge 2$, i.e.,~$e'$ joins vertices in~$V_i^3$ and~$V_j^3$, where~$2 \le i < j \le \ell -1$. In the former case, let us denote by~$e$ the edge of~$R$, to which~$e'$ is adjacent to, and observe that~$e$ is unique. This means that 3 vertices of the~$\mathcal{K}_\ell$'s that we want to count are specified (the one endpoint of~$e$, the intersection of~$e$ and~$e'$, and the other end of~$e'$); hence, the number of~$\mathcal{K}_\ell$'s with~$e$ and~$e'$ is at most~$n^{\ell-3}$. In the latter case we want to count~$\mathcal{K}_\ell$'s in~$G_{n,p}$, which have a vertex in~$V_i^3$ and~$V_j^3$, where~$2 \le i < j\le \ell-1$. These~$\mathcal{K}_\ell$'s have exactly one vertex in each~$V_x^3$, such that~$x\not\in\{1,i,j\}$, and two vertices in~$V_1^3$, which are endpoints of an edge in~$R$. As the number of admissible indexes~$x$ is~$\ell-4$, the number of~$\mathcal{K}_\ell$'s with~$e'$ is at most~$|R| \cdot n^{\ell-4} \le n^{\ell-3}$ (observe that trivially~$|R|\le n$, as~$R$ is a matching). 

By combining the last two results we conclude that the number of~$\mathcal{K}_\ell$'s, which contain exactly one edge from~$R$ and are otherwise edge-disjoint, is at least
\[
	\frac{|H_3| \cdot p^{3\ell^2-1}\cdot n^{\ell-2}}{\log n} \cdot \frac{1}{1 + (\binom{\ell}{2}-1)\cdot n^{\ell-3}} 
	\ge |H_3| \cdot \frac{p^{3\ell^2}}{\log n}n,
\]
for $n$ sufficiently large,  which is the desired lower bound for the number of missing edges.

To complete the proof we show (vi). Let~$d$ be the maximum horizontal degree of a vertex~$v\in X$. As the total number $|H|$ of horizontal edges is at most~$d\cdot |X| + |H_3|$. By exploiting (i)-(v) this implies
\begin{equation}\label{as:e1}
|H|\le d\cdot 2p^{-\ell^2-2}\log n  + p^{-5\ell^2}\log^2n,
\end{equation}
It thus clearly suffices to show that $d \le p^{-3\ell^2}\log n$. Suppose that
$d \ge cp^{-2\ell^2}\log n$ (where $c$ is the constant $c=c(\ell-1)$ from Proposition~\ref{prop:KellhelperII}). Without loss of generality we may assume that~$v\in V_1$. Then the degree of~$v$ in every~$V_i$ is also at least~$d$, as otherwise the chosen partition would have not been maximal. By applying Proposition~\ref{prop:KellhelperII} we readily obtain that the number of edge-disjoint~$\mathcal{K}_{\ell-1}$'s joining the neighborhoods of~$v$ is at least~$c'pd^2$, for an appropriate constant~$c'>0$, 
which implies that there are at least that many missing edges. As the number of missing edges is at most the number of horizontal edges, this together with $(\ref{as:e1})$ implies that
\[
	c'pd^2 \le 2p^{-\ell^2-2}\log nd + p^{-5\ell^2}\log^2n,
\]
from which we easily deduce~$d \le p^{-3\ell^2}\log n$, as desired. This completes the proof.
\end{proof}
{\bf Proof of Theorem~\ref{thm:cliques}} Theorem~\ref{thm:cliques} can be proved in a completely analogous way as Theorem~\ref{thm:extremalMain} (see Section~\ref{sec:FKGproof}). The definitions of several events, as well as the partial ordering of graphs with respect to partitions of the vertex set all generalize in an obvious and natural way from bipartitions to~$(\ell-1)$-partitions. We leave the straightforward details to the reader.

\bibliographystyle{alpha}
\footnotesize
\bibliography{ESRG}

\end{document}